\documentclass[12pt]{article}

\usepackage{graphicx}
\usepackage{hyperref}
\usepackage{amsmath,epsfig,epsf,psfrag,natbib}
\usepackage{amssymb}
\usepackage[latin1]{inputenc}
\usepackage{rotating}
\usepackage{subfigure}
\usepackage{amsbsy}
\usepackage{lscape,url,epstopdf}
\usepackage{color,latexsym,amsfonts,amsthm,amscd,graphicx}
\usepackage{mathrsfs}
\usepackage{float, multirow, booktabs, threeparttable, diagbox}

\newcommand{\ben}{\begin{enumerate}}
	\newcommand{\een}{\end{enumerate}}
\newcommand{\be}{\begin{equation}}
	\newcommand{\ee}{\end{equation}}
\newcommand{\bas}{\begin{eqnarray*}}
	\newcommand{\eas}{\end{eqnarray*}}
\newcommand{\ba}{\begin{eqnarray}}
	\newcommand{\ea}{\end{eqnarray}}
\newcommand{\bit}{\begin{itemize}}
	\newcommand{\eit}{\end{itemize}}

\newtheorem{coro}{Corollary}
\newtheorem{theorem}{Theorem}
\newtheorem{lemma}{Lemma}
\newtheorem{example}{Example}

\newtheorem{remark}{Remark}
\newtheorem{definition}{Definition}
\newtheorem{assumption}{Assumption}

\newcommand{\e}{ { \mathbb{E}}}

\newcommand{\betap}{ { \boldsymbol{\eta}_{P}}}
\newcommand{\betapt}{ { \boldsymbol{\eta}^{0}_{P}}}

\newcommand{\betaq}{ { \boldsymbol{\eta}_{Q}}}
\newcommand{\bPiQ}{\boldsymbol{\pi}_{Q}}

\newcommand{\hatbPiQ}{\widehat{\boldsymbol \pi}_{Q}}

\def\T{{ \mathrm{\scriptscriptstyle \top} }}

\newcommand{\bm}[1]{\mbox{\boldmath{$#1$}}}
\newcommand{\hatbetap}{\widehat{\boldsymbol \eta}_{P}}

\newcommand{\bphit}{{\boldsymbol \phi}^{0}}

\usepackage{cleveref}
\usepackage{verbatim}
\usepackage{hyperref}
\usepackage{authblk}
\usepackage{makecell, arydshln}
\newcommand{\bxi}{\boldsymbol{\xi}}

\usepackage{authblk}

\begin{document}
	\date{}
	\title{
	Unsupervised optimal deep transfer learning for classification under  general conditional shift
	}
\author[1]{Junjun Lang}
\author[1]{Yukun Liu\thanks{Corresponding author:  ykliu@sfs.ecnu.edu.cn}}
{\small   \affil[1]{ \small
		KLATASDS-MOE,  School of Statistics,
		East China Normal University,
		Shanghai 200062, China}
}
	\maketitle

\abstract{
%Classifiers trained on labelled source data may yield misleading results when applied to unlabelled target data from a different distribution. Transfer learning can rectify this by transferring knowledge from source to target data, but its validity often hinges on strict assumptions like label shift. In this paper, we introduce a general conditional shift (GCS) assumption, which is novel and encompasses label shift as a special case. Under GCS, we show that the target distribution and the shift function are identifiable.We propose to  estimate the conditional probabilities $\eta_P$ for source data by leveraging deep neural networks (DNNs).  Following the transfer of the DNN estimator, we  estimate the target label distribution ${\bm\pi}_Q$ using a pseudo-maximum likelihood technique. Ultimately, we construct our proposed Bayes classifier by incorporating these estimates and avoiding the need to estimate the shift function. Concentration bounds for our estimators of  both $\eta_P$  and ${\bm\pi}_Q$ are establish in terms the intrinsic dimension of $\eta_P$. We show that our DNN-based classifier achieves the optimal minimax rate up to a logarithm factor. A notable advantage of our method is its ability to effectively mitigate the curse of dimensionality when $\eta_P$  possesses a low-dimensional structure. Numerical simulations and the analysis of an Alzheimer's disease dataset showcase its outstanding performance.

Classifiers trained solely on labeled source data may yield misleading results when applied to unlabeled target data drawn from a different distribution. Transfer learning can rectify this by  transferring knowledge from source to target data, but its effectiveness frequently relies on stringent assumptions, such as label shift. In this paper, we introduce a novel General Conditional Shift (GCS) assumption, which encompasses label shift as a special scenario. Under GCS, we demonstrate that both the target distribution and the shift function are identifiable. To estimate the conditional probabilities ${\bm\eta}_P$ for source data, we propose leveraging deep neural networks (DNNs). Subsequent to transferring the DNN estimator, we estimate the target label distribution ${\bm\pi}_Q$ utilizing a pseudo-maximum likelihood approach. Ultimately, by incorporating these estimates and circumventing the need to estimate the shift function, we construct our proposed Bayes classifier. We establish concentration bounds for our estimators of both ${\bm\eta}_P$  and ${\bm\pi}_Q$ in terms of the intrinsic dimension of ${\bm\eta}_P$ . Notably, our DNN-based classifier achieves the optimal minimax rate, up to a logarithmic factor. A key advantage of our method is its capacity to effectively combat the curse of dimensionality when ${\bm\eta}_P$  exhibits a low-dimensional structure. Numerical simulations, along with an analysis of an Alzheimer's disease dataset, underscore its exceptional performance.
}


\section{Introduction}

Classifiers trained using labelled source data are typically employed to classify unlabelled target data that follow the same distribution. When the target data distribution deviates from that of the source data, the classifiers may lead to biased or misleading conclusions. Transfer learning is capable of correcting such classifiers by integrating information from both the labelled source data and the unlabeled target data. Generally speaking, transfer learning, also known as domain adaptation, arises in situations where the availability of data from the target domain is limited, whereas a considerable amount of data exists in a related source domain, thereby enabling knowledge transfer from the latter to the former. It has attracted significant attention across a wide variety of applications, particularly in computer vision \citep{tzeng2017adversarial}, speech recognition \citep{huang2013cross}, and machine learning \citep{iman2023review}. See \cite{weiss2016survey} and \cite{zhuang2020comprehensive} for more information regarding transfer learning and its applications. This paper focuses on transfer learning within the framework of nonparametric classification and is committed to utilizing modern deep learning techniques to develop an efficient algorithm that adapts classifiers trained on the source to the target domain under a general distribution shift assumption.

\subsection{Problem statement}

Let $X \in \mathcal{X}\subset\mathbb{R}^{d}$ for $d\geq 1$
 denote a vector of covariates or features and let $Y\in \mathcal{Y} \subset \mathbb{R}$
represent a discrete outcome or label.
Denote the joint distributions of $(X, Y)$ from  the source  and target domains
by   $P_{(X, Y)}$ and $Q_{(X, Y)}$, respectively.
Suppose that we have labelled source data and  unlabelled or unsupervised target data,
\ba
\label{data}
 (X_{1}^{P}, Y_{1}^{P}),\cdots,(X_{n_{P}}^{P}, Y_{n_{P}}^{P})\stackrel{i.i.d}{\sim}P_{(X,Y)},
 \qquad X_{1}^{Q},\cdots,X_{n_{Q}}^{Q}\stackrel{i.i.d}{\sim}Q_{X},
\ea
where $Q_{X}$ is the distribution function or probability measure of $X$ in the target domain.
Because the labels of the target data  are unavailable,
it is challenging, if not impossible, to  learn  a  classifier
with theoretical guarantee for the target data based solely on the  target data themselves.
This goal can be achieved by transferring information from the source data to the target data.
A typical example is medical diagnosis, where $X$ usually represents the physical indicators
of patients that are relatively easy to measure,
while $Y$ measures the quality  of the patient's condition after treatment,
which is often costly or difficult to measure.
To save costs, doctors may wish to infer whether patients can end treatment
based on the data of previous patients combined with the easily measurable physical indicators of the current patients.

 Transfer learning or domain adaptation is impossible without underlying assumptions \citep{Ben-David2010}.
Its validity  hinges on the similarity between the distributions of source and target data.
Within this framework, two prominent similarity assumptions
are covariate shift \citep{Kpotufe2021,Sugiyama2008} and label shift \citep{Saerens2002adjusting, maity2022minimax}.
In the scenario of covariate shift,  the conditional distributions  of $Y$ given $X$
 remain consistent between the source and target data.
 Consequently, classifiers trained using source data yield valid predictions when directly employed on target data.
% therefore   classifiers trained with the  source   data will produce valid results when directly applied to   the    target data.
Let $P_{X\mid Y}(x\mid y)$ and $Q_{X\mid Y}(x\mid y)$  be the probability density or mass functions of  $X$ given $Y=y$ in the source and target  domains, respectively. Conversely, label shift posits that
$
P_{X\mid Y}(x\mid y) = Q_{X\mid Y}(x\mid y)
$
and the marginal distributions of $Y$  may diverge between the source and target domains.
Label shift typically arises when external factors alter the distribution of labels.
For instance, in medical diagnosis, symptoms (features) are often indicative of diseases (labels),
but the distribution of diseases itself may shift.
%Label shift usually occurs when the distribution of labels changes due to certain external factors. For example, in medical diagnosis,  symptoms (features) is usually determined by diseases (labels), however the distribution of diseases may change.
%diseases (labels) causes symptoms (features).
%In disease classification and medical diagnosis, the proportion of lung cancer and stomach cancer cases in the original training dataset is comparable. However,  due to a decrease in the number of smokers and changes in dietary habits, there has been a decrease in lung cancer cases and an increase in stomach cancer cases in the test population, which leads to label shift.
%In the example of medical diagnosis,    the proportion of good quality of the patient's condition after treatment may be greater than half before the  COVID-19 pandemic. As COVID-19 may reduce people's immunity, there may be a decrease in this proportion after  the  COVID-19 pandemic, resulting in label shift.

Despite its widespread popularity in the realm of transfer learning, the label shift assumption may be restrictive, prohibiting any alterations in the conditional distribution of $X$ given $Y$. For instance, the advent of COVID-19 has demonstrated that, even when a patient's disease is known, the distribution of their symptoms can change post-pandemic. To accommodate such variations, we propose relaxing the label shift assumption and adopting a    more general distribution shift assumption:
\begin{assumption}[General Conditional Shift, GCS] \label{assumption-GCS}
There exists a function class $\mathscr{H}$  on $\mathcal{X}$
  containing
  the constant function 1 such that for a function   $h(\cdot) \in \mathscr{H}$,
\ba
\label{GCS}
\frac{Q_{X\mid Y}(x\mid y)}{P_{X\mid Y}(x\mid y)} =  h(x), \quad \forall x\in\mathcal{X}, \;  y\in\mathcal{Y}.
\ea
\end{assumption}
\noindent
The GCS assumption stipulates that the ratios of the conditional densities of $X$ given $Y=y$ in the target and source domains are equivalent to a common function $h(x)$ for all values of $y$. When the function class $\mathscr{H}$ comprises solely the element 1, this assumption simplifies to label shift. Furthermore, it can be perceived as a modification of the conditional shift proposed by \cite{liu2021adversarial} or the exponential tilting model presented by \cite{maity2023understanding}.

\subsection{Our contributions}

In this paper, we develop a minimax optimal classifier for unsupervised target data under the GCS  assumption by leveraging knowledge transferred from the source domain using deep neural networks (DNNs). Our contributions are outlined as follows:

\ben

\item
We introduce  a novel distribution shift assumption, namely the GCS model presented in \eqref{GCS}, which encompasses the traditional label shift assumption as a special case. The GCS model not only provides greater flexibility compared to label shift but also enables the identification and consistent estimation for the classification function of target domain by utilizing information from both source and target data.

\item
Under the GCS assumption, we incorporate DNNs into transfer learning for classification tasks based on labeled source data and unlabeled target data. In contrast, traditional non-parametric methods, such as kernel-based techniques \citep{maity2022minimax} or the K-NN method \citep{TonyCai2021}, have been utilized for this purpose in the existing literature. Given the DNN estimator of $\betap(x)$, the classification function of source domain,  we introduce a pseudo-maximum likelihood estimator $\hat{\bm \pi}_Q$ for the marginal distribution or class proportion ${\bm \pi}_Q$ of the label $Y$ in the target domain. We establish concentration bounds for both the DNN estimator of $\betap(x)$ and $\hat{\bm \pi}_Q$. Subsequently, the proposed classifier, which serves as an estimator of the Bayes classifier for the target data, is constructed by integrating these two estimators. One of the most fascinating features of our method lies in its ingenious way of circumventing the necessity of estimating the unknown function $h(\cdot)$ in the GCS model \eqref{GCS}.

\item
We assume that the regression or classification function of the source domain belongs to the {\it composite smoothness function class} \citep{Schmidt-Hieber2020}, and demonstrate that the convergence rate of our proposed classifier is contingent upon the intrinsic dimension of the ideal classification function. Our method is capable of mitigating the curse-of-dimensionality when the classification function in the source domain exhibits a certain low-dimensional structure \citep{Schmidt-Hieber2020, Zhong2022}. In contrast, existing works on transfer learning for classification \citep{maity2022minimax, TonyCai2021, liu2020computationally} typically grapple with the curse-of-dimensionality issue, where the convergence rates of their classifiers can be extremely slow in the presence of high-dimensional covariates.

\item
We derive the minimax optimal convergence rate for classifiers in the target domain in terms of excess risk in general multi-class settings. Our findings reveal that the proposed DNN-based classifier attains the optimal minimax rate, asymptotically up to a logarithmic factor. These results constitute novel contributions to the literature on unsupervised transfer learning for classification, particularly in scenarios involving labeled source data and unlabeled target data.

\een

In summary, under the GCS model and within a multi-class setting, we have established theoretical guarantees for our DNN-based estimators concerning the classification function in the source domain and class proportions in the unlabeled target domain. Leveraging the strengths of deep transfer learning, we have developed an efficient classifier specifically designed for the target domain. This classifier has demonstrated its prowess in addressing the curse of dimensionality, achieving minimax optimality, and exhibiting remarkable efficacy in a real-world application.

\subsection{Related works }

The GCS assumption is intimately tied to the concept of label shift. Within the context of unsupervised transfer learning for classification under label shift, the estimation of the label distribution (or class proportions) of the target data, as well as the distribution shift function between the target and source domains, represents a pivotal task. Broadly speaking, existing methodologies for this purpose can be categorized into two primary groups.
The first category revolves around comparing the moments of the target data with equivalent expressions derived from the source data distribution and the label distribution of the target data. Notable examples include the maximum mean discrepancy-based method \citep{zhang2013domain, iyer2014maximum}, the black box shift learning (BBSL) method, and the regularized learning under label shift (RLLS) methods \citep{lipton2018detecting, azizzadenesheli2019regularized}.
The second category encompasses the maximum likelihood label shift (MLLS) method \citep{Saerens2002adjusting, alexandari2020bias}, which is an instance of the expectation-maximization (EM) algorithm \citep{Dempster1977maximum}. \cite{garg2020unified} presented a unified framework that connects MLLS and BBSL, revealing that MLLS often surpasses BBSL and RLLS in performance. Additionally, \cite{jain2016estimating} investigated the same problem, albeit with source data that can be considered unlabeled.

Another crucial task is the construction of valid classifiers for the target domain, as explored in \citep{Saerens2002adjusting, azizzadenesheli2019regularized, maity2022minimax}. Two predominant approaches include weighted empirical risk minimization \citep{lipton2018detecting, azizzadenesheli2019regularized} and the plugging-in method (in conjunction with the Bayes classifier) \citep{maity2022minimax}. Many existing studies presuppose the availability of a classification model or the conditional distribution of $Y$ given $X$ for the source data, which can be estimated using an artificial neural network, logistic regression, or kernel-based methods. However, these studies either lack theoretical guarantees \citep{Saerens2002adjusting, alexandari2020bias} or may exhibit suboptimal finite-sample performance \citep{maity2022minimax}.
Recently, under the label shift assumption, \cite{Lee2024doubly} developed a doubly flexible estimation method for the target population within the framework of semiparametric efficiency theory. Nonetheless, their focus was on parameter estimation rather than classification. Most importantly, the existing results derived under the label shift assumption may not be applicable under the GCS assumption.

Another closely related distribution shift assumption to GCS is posterior drift, which signifies that the conditional probabilities of labels given features differ between the source and target domains. Existing research on transfer learning under posterior drift has been conducted under specific similarity assumptions \citep{TonyCai2021, reeve2021adaptive, maity2024linear}. Our GCS bears a strong resemblance to the linear adjustment model proposed by \cite{maity2024linear}, which represents a particular similarity assumption between the source and target domains. However, unlike the current paper, where the observed target data are unlabeled, these prior works necessitate access to target labels and do not grapple with identifiability challenges. Consequently, they are not applicable to our specific setting.

 A pivotal component of our methodology is the deep neural network (DNN), which has driven substantial advancements over the past few decades \citep{szegedy2015going, sarikaya2014application, miao2018coronary}. Mathematically, a DNN is a composite function comprising simplistic functions layered in multiple tiers, imparting it with the remarkable ability to approximate complex relationships with high precision. Theoretical analyses have demonstrated its capacity to approximate a wide spectrum of functions with specific structures, such as composite smooth functions \citep{BauerandKohler2019}, functions in Sobolev spaces \citep{Yarotsky2017}, and non-smooth functions \citep{ImaizumiandFukumizu2019}.
Owing to its formidable approximation power, DNNs have garnered significant attention in numerous fields, including nonparametric regression \citep{Schmidt-Hieber2020, jiao2023deep}, survival analysis \citep{Zhong2022}, quantile regression \citep{shen2024nonparametric}, and causal inference \citep{chen2024causal}, among others. In this paper, we utilize DNNs to transfer knowledge between domains in an unsupervised setting under the GCS model, representing a variant of deep transfer learning \citep{tan2018survey, iman2023review}.

\subsection{Organization}

The organization of the remainder of this paper is as follows. Section \ref{sec-notations-identifiability} establishes the identifiability conditions for the target distribution within the framework of our GCS model and formally introduces a nonparametric setting for classification.
In Section \ref{sec-method}, we present the proposed DNN estimator for the classification function of the source data, along with a pseudo-maximum likelihood estimator for the class proportions of the target data. By integrating these two estimators, we derive the proposed classifier for the target data.
Section \ref{sec-theory} undertakes a comprehensive examination of the large-sample properties of the proposed estimators and classifier, encompassing concentration bounds, upper bounds in terms of excess risk, and minimax optimality.
Sections \ref{sec-simu} and \ref{sec-real} provide a detailed simulation study and a real data analysis, respectively.
Finally, Section \ref{sec-con} concludes with some remarks. For clarity and ease of reading, all technical proofs have been deferred to the supplementary material.

\section{Identifiability  and nonparametric  classification}
\label{sec-notations-identifiability}

Without loss of generality, we assume the label space
$\mathcal{Y}$  is given by $ \{1,2,\cdots, k\}$ for $k\geq2$,
and the feature space $\mathcal{X}$ is $[0,1]^{d}$.
Our objective is to develop a valid classifier   using a deep neural network (DNN) for the unsupervised target data,
based on the observed  data  specified in \eqref{data},
 under the GCS assumption  in \eqref{GCS}.

\subsection{Notation}

In this paper, we adopt  a nonparametric assumption  concerning the   classification of the target data.
To facilitate our discussion,  we introduce the necessary  notation.
Let $P$ and $Q$ represent the probability measures induced by $P_{(X, Y)}$ and $Q_{(X, Y)}$, respectively,
and let  $P_{X}$ and $Q_{X}$ denote the distribution functions or probability measures
of $X$ in   the source and target domains, respectively.
We  use $\mathcal{P}$ to denote the pair of   source and target distributions  $(P_{(X,Y)}, Q_{(X,Y)})$
and define $[k] = \{1, 2, \ldots, k\}$.
For $  l \in [k]$, let  $\pi_{Q,l} = Q(Y=l)$, $\pi_{P,l} = P(Y=l)$,
$\eta_{Q,l}(x) = Q(Y=l \mid X=x)$, $ \eta_{P,l}(x) =  P(Y=l \mid X=x)$,
$f_{l}(x)  =  Q_{X\mid Y}(x\mid l)$ and $\alpha_{l}  =   \log(\pi_{P,k}/\pi_{P,l})$.
Additionally, we  define
$\betap(x)  =  (\eta_{P,1}(x),\cdots,\eta_{P,k-1}(x))^{\top}$,
$\betaq(x)  =   (\eta_{Q,1}(x),\cdots,\eta_{Q,k-1}(x))^{\top}$ and ${\boldsymbol \pi_{Q}}  =  (\pi_{Q,1},\cdots, \pi_{Q,k-1})^{\top}$.

We use  $\mathcal{D} = \{ (X_{1}^{P}, Y_{1}^{P}),
\cdots,(X_{n_{P}}^{P}, Y_{n_{P}}^{P}); X_{1}^{Q},\cdots,X_{n_{Q}}^{Q}\}$
to denote the observed  data   \eqref{data} in
$\mathcal{S} = (\mathcal{X}\times\mathcal{Y})^{n_{P}}\times \mathcal{X}^{n_{Q}}$,
and use    $\mathcal{H} = \{\tilde{h}:\mathcal{X}\mapsto\mathcal{Y}\}$
to denote the set of all classifiers on $\mathcal{X}$.
For any classifier $\tilde{h}\in\mathcal{H}$, its excess-risk  is defined as
\ba
\label{definition-excess-risk}
\mathcal{E}_Q(\tilde{h})=Q(Y \neq \tilde{h}(X))-Q\left(Y \neq f_Q^*(X)\right),
\ea
where $f_Q^*(X) = \operatorname*{argmax}_{l\in\mathcal{Y}}\eta_{Q, l}(X)$ is the  Bayes
classifier for the target data.

Let  $\|\cdot\|$ denote the standard Euclidean norm.
We employ $\lambda_{\text{min}}(\cdot)$ to
represent the minimal eigenvalue  of a matrix, and
use $\|\cdot\|_{\infty}$ to signify the supremum norm of  matrices,  vectors or  functions.
Furthermore, let $\mathbb{N}$, $\mathbb{N}_{+}$ and $\mathbb{R}_{+}$
denote the sets of natural numbers, positive natural numbers, and positive real numbers, respectively.
For any distribution or probability measure  $M$
and any function $f$ defined on $\mathcal{X}$,
$\mathbb{E}_{M}[f(X)]$ denotes the expectation value of $f(X)$
with respect to $X\sim M$.
For any values $a$ and $b$, let $a\vee b := \max\{a, b\}$ and
$a\wedge b := \min\{a, b\}$,  where $:=$ indicates ``defined as''.

\subsection{Identifiability}

 Because   the target data lack labels,
the target distribution $Q_{(X,Y)}$ and  $h(\cdot)$
in \eqref{GCS}  may not be identifiable without incorporating  additional assumptions.
 Therefore, it is crucial to examine the issue of identifiability prior to presenting our learning algorithm.

\begin{assumption}
	\label{assumption-identifiability} (i)   For the function class $\mathscr{H}$ in Assumption \ref{assumption-GCS}
 and any $\tilde{h}\in\mathscr{H}$,
  $ \tilde{h}(x)\cdot P_{X\mid Y}(x\mid 1), \ldots, \tilde{h}(x)\cdot P_{X\mid Y}(x\mid k)$  are nondegenerate probability density or mass functions. (ii) As functions of $x$ on $\mathcal{X}$,  $P_{X\mid Y}(x\mid 1), \ldots,  P_{X\mid Y}(x\mid k)$ are linearly independent. (iii) $\pi_{Q, i} >0$  for all $1\leq i\leq k$.
\end{assumption}

Assumption \ref{assumption-identifiability}(i)-(ii) are relatively straightforward.
Under Assumption \ref{assumption-GCS}, the constant function $1$  belongs to the function class $\mathscr{H}$.
Then,   Assumption \ref{assumption-identifiability}(i) implicitly
assumes that   the conditional probability density or mass functions $P_{X\mid Y}(x\mid 1), \ldots, P_{X\mid Y}(x\mid k)$
are all nondegenerate.
 In Assumption \ref{assumption-identifiability} (iii),
we impose the condition that  $\pi_{Q, i} >0$ for all $1\leq i\leq k$
 to ensure that the target distribution has the same support as the source distribution.

We have discovered that the combination of Assumptions \ref{assumption-GCS} and  \ref{assumption-identifiability} constitutes a sufficient condition for the identifiability of $Q_{(X, Y)}$ and   $h(\cdot)$.
This finding  is formally articulated in the following lemma.

\begin{lemma}
	\label{prop-identifiability-GCS}
Under Assumptions \ref{assumption-GCS} and  \ref{assumption-identifiability},
the target distribution $Q_{(X, Y)}$
and the function $h(\cdot)$ in model \eqref{GCS} are identifiable
from the data in \eqref{data}.
\end{lemma}

%In what follows, for ease of discussion,  we directly assume that the
%$Q_{(X, Y)}$ and   $(\alpha, \beta)$ is identifiable. Otherwise,  this assumption
%can be simply replaced with   Assumption \ref{assumption-identifiability}.

\subsection{Nonparametric setting for classification}

\label{sec-setting}

Before presenting our nonparametric setting for classification, we introduce two types of function classes that play crucial roles in the construction of our classifier. The first type is H$\rm \ddot{o}$lder classes, which are frequently utilized in nonparametric regression and machine learning \citep{stone1982optimal, TonyCai2021, maity2022minimax, reeve2021adaptive}.

\begin{definition}[H$\rm \ddot{o}$lder classes]
Let  $a$ and $M$ be positive constants, and let $\mathbb{D}\subset\mathbb{R}^{d}$. The following  class
\bas
\mathcal{H}_{d}^{a}(\mathbb{D}, M)
=\left\{f:\mathbb{D}\mapsto\mathbb{R}\bigg|\sum\limits_{\lambda:|\lambda|<a}\|\partial^{\lambda}f\|_{\infty}
+\sum\limits_{\lambda:|\lambda|=\lfloor a \rfloor}\sup_{x,y\in\mathbb{D},x\neq y}\frac{|\partial^{\lambda}f(x)-\partial^{\lambda}f(y)|}{\|x-y\|_{\infty}^{a-\lfloor{a}\rfloor}}\leq M\right\}
\eas
is called a H$\rm \ddot{o}$lder class of smooth functions.
Here,
$\lambda=(\lambda_{1},\cdots,\lambda_{d})$ with  $\lambda_i$ being nonnegative integers,
$|\lambda|=\sum_{i=1}^{d}\lambda_{i}$, $\partial^{\lambda} := \partial^{\lambda_{1}}\cdots\partial^{\lambda_{d}}$, and
$\lfloor a \rfloor$ denotes the largest integer strictly smaller than $a$.
\end{definition}

A function within the class  $\mathcal{H}_{d}^{a}(\mathbb{D}, M)$
possesses a H$\rm \ddot{o}$lder smoothness index of  $ a$.
If a conditional mean regression function belongs to
 $\mathcal{H}_{d}^{a}(\mathbb{D}, M)$,
 \cite{stone1982optimal} revealed that the optimal minimax rate of convergence
for estimating it is $n_P^{-\frac{a}{2a+d}}$.
In numerous  transfer learning studies,
specific nonparametric components of the underlying problem  are presumed to belong  to
a  H$\rm \ddot{o}$lder class, such as  $\mathcal{H}_{d}^{a}(\mathbb{D}, M)$
\citep{TonyCai2021, maity2022minimax, reeve2021adaptive}.
Under this assumption, the  minimax optimal rate of convergence for the  classifier on the target data, in terms of excess-risk,  hinges on
 the smoothness index  $a$ and the dimension $d$.
 For instance, this rate may be
   $ n_P^{   - \frac{ a  (1+c_2)}{ 2 c_1 a +d }} + n_Q^{- \frac{a (1+c_2)}{  2a +d } }$
for some positive constants $c_1$ and $c_2$ \citep{TonyCai2021}.
When  the dimension $d$ is substantial, this rate  can be extremely slow.
To mitigate this curse of dimensionality,
we consider another type of  function class.

\begin{definition}[Composite smoothness function class  \citep{Schmidt-Hieber2020}]
	Let $q\in\mathbb{N}$, $M\geq1$, $\theta=(\theta_{0},\cdots,\theta_{q})\in\mathbb{R}_{+}^{q+1}$ and ${\bf k}=(k_{0},\cdots, k_{q+1})\in\mathbb{N}_{+}^{q+2}$,and
 ${\bf \tilde{k}}=(\tilde{k}_{0},\cdots,\tilde{k}_{q})\in\mathbb{N}_{+}^{q+1}$ with $\tilde{k}_{j}\leq k_{j}$, $j=0,\cdots, q$.
The following class
\bas
\mathcal{G}(q,\theta,{\bf k}, {\bf \tilde{k}}, M)
&= & \left\{f=g_q \circ \cdots \circ g_0\bigg| g_i=\left(g_{i j}\right)_j:\left[a_i, b_i\right]^{k_{i}} \rightarrow\left[a_{i+1}, b_{i+1}\right]^{k_{i+1}},\right. \\
			&& \hspace{2cm} \left.g_{i j} \in \mathcal{H}_{\tilde{k}_{i}}^{\theta_i}\left(\left[a_i, b_i\right]^{\tilde{k}_{i}}, M\right), \text { for some }\left|a_i\right|,\left|b_i\right| \leq M\right\}
\eas
with  $k_{0}=d$, $k_{q+1}=1$, $a_{0}=0$, $b_0=1$, and $f\circ g  = f(g(\cdot))$ is called
a composite smoothness function class.
\end{definition}

For a  constant vector   $\theta$,  let
$\gamma_{n_{P}} = \max_{i=0,\cdots, q}n_{P}^{-\tilde{\theta}_{i}/(2\tilde{\theta}_{i}+\tilde{k}_{i})}$,
where  $\tilde{\theta}_{i} = \theta_{i}\prod_{j=i+1}^{q}(\theta_{j}\wedge1)$. Let $i^{\star}\in\operatorname*{argmin}_{i\in \{0,\cdots, q\}}\tilde{\theta}_{i}/(2\tilde{\theta}_{i}+\tilde{k}_{i})$. To some extent,  $\tilde{k}_{i^{\star}}$ characterizes the intrinsic dimension
of the function class $\mathcal{G}(q,\theta,{\bf k}, {\bf \tilde{k}}, M)$.
In this paper, we assume that after an equivalence transformation,
the ideal classification function of the source domain
belongs to a composite smoothness function class.
We uncover that under this assumption, the convergence rate of our estimator (See Section \ref{Sec-DNN-estimator}) for the classification function  hinges upon the intrinsic dimension of the classification function.  Intuitively,  the structure of the composite functions within $\mathcal{G}(q,\theta,{\bf k}, {\bf \tilde{k}}, M)$ bears resemblance to
 the multi-layer perceptron structure of neural networks, which endows the resultant estimator with relatively smaller approximation errors or a faster convergence rate.
Consequently, we anticipate that the proposed methodology can mitigate the curse of dimensionality when the classification function in the source domain exhibits a low-dimensional structure \citep{Schmidt-Hieber2020,Zhong2022}.

With the aforementioned preparations in place, we impose the requirement that  the distribution pair $\mathcal{P}$ of interest must adhere to   Assumption \ref{assume-nonparametric-setting}.
 \begin{assumption}
 	\label{assume-nonparametric-setting}
 The distribution pair $\mathcal{P}$ satisfies
  (i) Assumption \ref{assumption-GCS} (the GCS assumption) holds to be true.
 (ii)
 There exists a group of parameters  $(q,\theta,{\bf k}, {\bf \tilde{k}}, M)$
 such that  $\phi_{k}(\cdot)  =  0$  and
 \ba
\label{model-regression-multi}
\eta_{P, l}(\cdot) = \frac{\exp\{\phi_{l}(\cdot)\}}{1+\sum\limits_{i=1}^{k-1}\exp\{\phi_{i}(\cdot)\}},~~\phi_{l}(\cdot)\in	\mathcal{G}(q,\theta,{\bf k}, {\bf \tilde{k}}, M), \quad   l\in [k-1].
\ea
 (iii) There exists $\tilde c>1$ such that $ dQ_{X}(x) / dP_{X}(x) \leq \tilde{c}$ for all $x\in\mathcal{X}$.
  (iv) Assumption \ref{assumption-identifiability} is satisfied.
 \end{assumption}

  Assumption \ref{assume-nonparametric-setting} (iii) is a weaker version of
   the strong  density assumptions  in Definition 4  of \cite{TonyCai2021}.
Given a set of  parameters $(h, q,\theta,{\bf k}, {\bf \tilde{k}}, M, \tilde{c})$,
we denote   the class of distribution pairs  $\mathcal{P}$ that satisfy
Assumption \ref{assume-nonparametric-setting}
as  $\mathcal{P}^{\star}(h, q,\theta,{\bf k}, {\bf \tilde{k}}, M, \tilde{c})$.
For brevity, we may also refer to this class as  $\mathcal{P}^{\star}$.
To signify the true value of parameter,
we append a suffix   $0$ to it; for instance,
 $\bPiQ^{0}$ represents the true value of $\bPiQ$.

\section{Classification with deep transfer learning}
\label{sec-method}
This section introduces the proposed nonparametric classifier that leverages
deep transfer learning.
Initially, we train  a {\it Deep Neural Network} (DNN) estimator, denoted as  $\hat{\bm \eta}_P$, for ${\bm \eta}_P$. Subsequently, we estimate   ${\bm\pi}_{Q}$ using
a pseudo maximum likelihood estimator, named $\hat{\bm\pi}_{Q}$.
The  classifier is constructed by combining
 $\hat{\bm \eta}_P$ and $\hat{\bm\pi}_{Q}$ through Bayes' formula.
A notable advantage  of our method is its innovative circumvention of the
necessity  to estimate
 the unknown function $h(\cdot)$ in the GCS model \eqref{GCS}.

\subsection{DNN estimator for ${\bm \eta}_P$ \label{Sec-DNN-estimator}}

Since ${\bm \eta}_P$ represents the regression of  $Y$ on $X$
under the source data distribution and the source data are fully available,
${\bm \eta}_P$ is nonparametrically identifiable
and can be consistently estimated using any nonparametric regression technique.
Within the framework of model \eqref{model-regression-multi},
estimating ${\bm \eta}_P$
is equivalent to estimating
  ${\boldsymbol \phi}(\cdot) =  (\phi_{1}(\cdot),\cdots, \phi_{k-1}(\cdot))^{\top}$.
  Given its remarkable flexibility and widespread popularity,
we propose  using a DNN to estimate   ${\boldsymbol \phi}(\cdot)$ and
consequently ${\bm \eta}_P$.

\begin{figure}[htbp]
	\centering
	\includegraphics[width=0.4\linewidth]{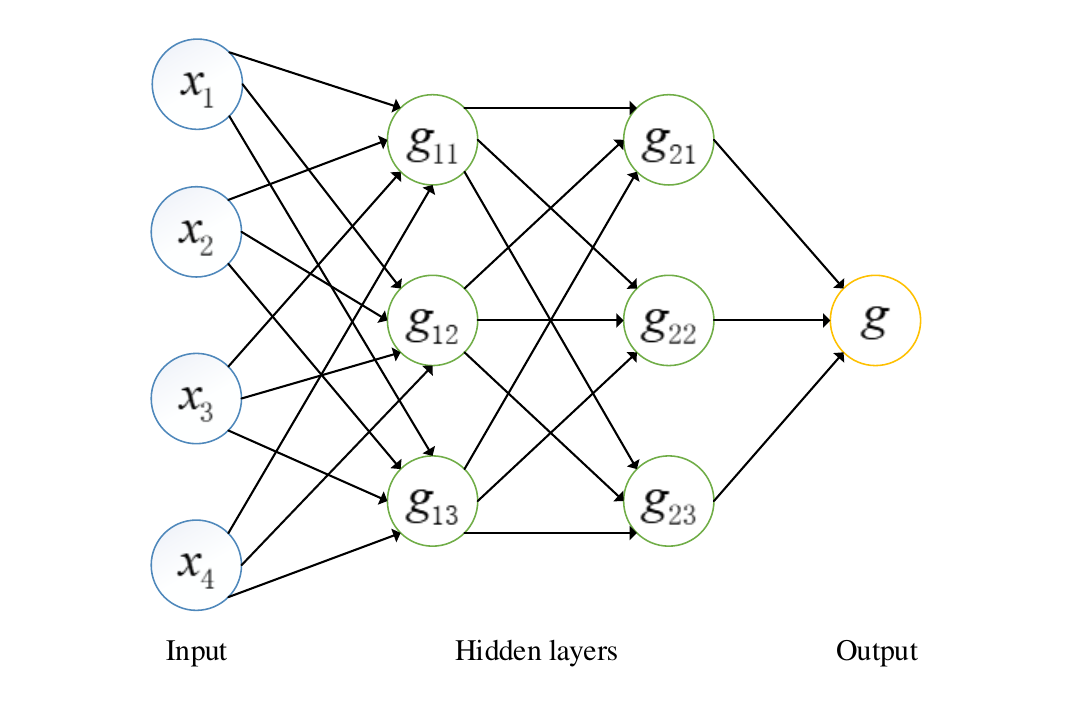}
	\label{figure-DNN}
	\caption{A three layers deep neural network with $K = 2$ and ${\bf p} = (4,3,3,1)$. }
\end{figure}

Let $K\in\mathbb{N}_{+}$  and ${\bf p}=(p_{0},\cdots, p_{K}, p_{K+1})\in\mathbb{N}_{+}^{K+2}$.
A $(K+1)$-layer DNN with layer width specified by ${\bf p}$ is a composition function
$g:\mathbb{R}^{p_{0}}\mapsto\mathbb{R}^{p_{K+1}}$  defined recursively as follows:
\ba
		\label{model-DNN}
		 g(x)=W_{K}g_{K}(x),\;
		 g_{K}(x)=\sigma\left(W_{K-1}g_{K-1}(x)+\mu_{K-1}\right),\cdots, g_{1}(x)=\sigma\left(W_{0}x+\mu_{0}\right).
	\ea
Here, the matrices $W_{j}\in\mathbb{R}^{p_{j+1}\times p_{j}}$ for $j=0,\cdots, K$,
 and the vectors $\mu_{l}\in\mathbb{R}^{p_{l+1}}$ for  $l=0,\cdots, K-1$,
constitute the parameters of the DNN.
The activation function  $\sigma$  operates component-wise,
meaning that  $g_{j}=(g_{j1},\cdots, g_{jp_{j}})^{\top}$
  maps from $\mathbb{R}^{p_{j-1}}$ to $\mathbb{R}^{p_{j}}$
for each $j=1,\cdots, K$.
The parameter $K$ represents the depth of the DNN,
while  ${\bf p}$ enumerates the width of each layer ($p_{0}$ corresponds to the dimension
of the input variable, $p_{1},\cdots, p_{K}$ denote the dimensions of
the $K$ hidden layers, and $p_{K+1}$ specifies the dimension of the output layer).
For illustrative purpose, Figure \ref{figure-DNN} depicts a three-layers DNN
with $K = 2$ and ${\bf p} = (4,3,3,1)$.

The activation function  $\sigma$ must be selected beforehand.
While sigmoidal  activation functions   were commonly used in  neural networks,
 the nonsigmoidal  rectified linear unit (ReLU), denoted as
$
 x \vee 0
$,
offers  computationally   advantages
for DNN \citep{Schmidt-Hieber2020}.
In terms of statistical performance,
the ReLU has been found to outperform  sigmoidal activation functions for
classification problems \citep{Glorot2011deep,Pedamonti2018comparison}.
In this paper,  we adopt the ReLU activation function and
 consistently set $p_{0}=d$ and $p_{K+1}=1$.

A fully-connected, deep feedforward network often suffers from overfitting.
To address this challenge, we consider a sparsely connected network approach
\citep{Hanetal2015, Srinivasetal2017, Schmidt-Hieber2020, Zhong2022}.
Given $K, s\in \mathbb{N}_{+}$, ${\bf p}\in\mathbb{N}_{+}^{K+2}$ and $D>0$, we define
a class of sparse neural networks as follows:
\bas
\mathcal{F}(K,s,{\bf p}, D)&=& \bigg\{ g \bigg| g~\text{is a}~(K+1)\text{-layer DNN}~\text{with width}~{\bf p}~\text{such that } \|g\|_{\infty}\leq D,\\
&&\hspace{1cm}
 \sum\limits_{i=0}^{K} (\|W_{i}\|_{0}+\|\mu_{i}\|_{0})\leq s,  \;\;\;
 \max_{j\in \{0, 1, \cdots, K\}}\{\|W_{j}\|_{\infty},\|\mu_{j}\|_{\infty}\}\leq 1 \bigg\},	
\eas
where $\mu_{K}$ is a $p_{K+1}$-dimensional zero
and $\|\cdot\|_{0}$ denotes the number of nonzeros entries of a matrix or vector.
We model ${\boldsymbol \phi}  =  (\phi_{1} ,\cdots, \phi_{k-1} )^{\top}$
using  a sparsely connected network, specifically,
${\boldsymbol \phi}\in	\{\mathcal{F}(K,s,{\bf p}, D)\}^{k-1}$.

Under the model specified in  \eqref{model-regression-multi}, we propose
to estimate ${\boldsymbol \phi} $ using
 the following nonparametric maximum likelihood estimator
\ba
\label{definie-NMPLE-multi}
\widehat{\boldsymbol \phi} = \operatorname*{argmax}_{{\boldsymbol \phi}\in	\{\mathcal{F}(K,s,{\bf p}, D)\}^{k-1} }\sum\limits_{i=1}^{n_{P}}\left\{
\sum_{l=1}^{k-1} I(Y_{i}^{P} = l) \phi_{l}(X_{i}^{P})-\log\left(1+\sum\limits_{j=1}^{k-1}\exp\{\phi_{j}(X^{P}_{i})\}\right)\right\},
\ea
where  $\widehat{\boldsymbol \phi}(\cdot)
= (\hat{\phi}_{1}(\cdot),\cdots, \hat{\phi}_{k-1}(\cdot))^{\top}$
and $  I(A) $ denotes the indicator function of an event $A$.
Once  $\widehat{\boldsymbol \phi} $ is obtained,
we can construct a DNN estimator for $\betapt$ by simply plugging in
the estimate, i.e.
$\hatbetap(\cdot) : = (\hat{\eta}_{P,1}(\cdot),\cdots,
\hat{\eta}_{P,k-1}(\cdot))^{\top}$,
where
\[
\hat{\eta}_{P,l}(x)=  \frac{\exp\{\hat{\phi}_{l}(x)\}}{  1+\sum_{j=1}^{k-1}\exp\{\hat{\phi}_{j}(x)\}},
\quad l \in [k-1], \quad x\in\mathcal{X}.
\]
%{\blue The depth $K$ and the width vector ${\bf p}$ control the complexity and flexibility of the neural network $\mathcal{F}(K,s,{\bf p}, D)$ and, consequently, the conditional probabilities  $\hatbetap(x)$.}
The depth $K$ and the width vector ${\bf p}$ govern  both the complexity and flexibility of the neural networks in $\mathcal{F}(K,s,{\bf p}, D)$, as well as those inherent in the conditional probabilities $\hatbetap(x)$.
These hyperparameters can be selected using cross-validation or held-out validation
in practical applications.

\subsection{Pseudo maximum likelihood estimator of  ${\bm\pi}_{Q}$ }
\label{sec-pmle}

In the second step, we propose to estimate  ${\bm\pi}_{Q}$
 by a pseudo maximum likelihood estimation method.
For a generic unlabeled  observation $X^{Q}$,
it follows a $k$-component mixture model  given by
$
\sum_{j=1}^{k}\pi_{Q,j}f_{j}(x).
$
Under the GCS model \eqref{GCS},   for $l \in [k-1]$, we have
\ba
\frac{f_{l}(x)}{f_{k}(x)} = \frac{  Q_{X\mid Y}(x\mid l) }{ Q_{X\mid Y}(x\mid k)}
&=& \frac{  P_{X\mid Y}(x\mid l)\cdot h(x) }{ P_{X\mid Y}(x\mid k)\cdot h(x)}\nonumber  \\
&=& \frac{  P(Y=l|X=x) /P(Y=l) }{ P(Y=k|X=x) /P(Y=k)}
= \frac{  \eta_{P,l} (x)  }{ \eta_{P,k} (x) }  \frac{\pi_{P,k} }{\pi_{P,l} }. \label{density-ratio-pmle}
\ea
Given model \eqref{model-regression-multi},
since  $\phi_k(x) \equiv0$, we obtain
$
f_l(x)
=    \exp\{ \alpha_{l} + \phi_{l}(x) \}f_k(x)
$
with $\alpha_{l} = \log (\pi_{P,k}/\pi_{P,l} )$,
for $l \in [k]$.
Therefore, the density of $X^{Q}$ can be expressed as
\bas
 f_k(x)  \sum_{j=1}^{k}\pi_{Q,j} \exp\{ \alpha_{j} + \phi_{j}(x) \}
 =
  f_k(x)  \bigg[
  \sum_{j=1}^{k-1}\pi_{Q,j} \exp\{ \alpha_{j} + \phi_{j}(x) \} + \bigg(1-\sum_{l=1}^{k-1}\pi_{Q, l}\bigg) \bigg].
\eas
This implies that
the  log-likelihood contribution of $X^{Q}$    is
$
 \log   f_k(X^{Q}) + \ell_{Q}(X^{Q};{\boldsymbol\alpha}, {\boldsymbol \phi}, {\boldsymbol \pi_{Q}}),
$
where  ${\boldsymbol \alpha} =  (\alpha_{1},\cdots, \alpha_{k-1})^{\top}$ and
\bas
\ell_{Q}(X^{Q};{\boldsymbol\alpha}, {\boldsymbol \phi}, {\boldsymbol \pi_{Q}}) =\log\left[\sum\limits_{l=1}^{k-1}\pi_{Q,l}\exp(\alpha_{l}+\phi_{l}(X^{Q}))+\left(1-\sum_{j=1}^{k-1}\pi_{Q, j}\right)\right].
\eas
Since the term $\log   f_k(X^{Q})$ does not  involve  $\pi_{Q, l}$,
we can  ignore it  and consider
 $\ell_{Q}(X^{Q};{\boldsymbol\alpha}, {\boldsymbol \phi}, {\boldsymbol \pi_{Q}}) $
 as
the  likelihood contribution of $X^{Q}$.
%{\red Alternatively, we may  construct an estimate   $\hat  f_k(X^{Q})$ of
%$f_k(X^{Q})$ based on the source data  and proceed with $\hat  f_k(X^{Q})$
%in place of $f_k(X^{Q})$.
%Again, the term  $\log  \hat f_k(X^{Q})$  can be ignored
%as it does not depend on the unknown parameters  ${\bm\pi}_{Q}$.}

The likelihood function,
$ \sum_{i=1}^{n_Q}\ell_{Q}(X_i^{Q};{\boldsymbol\alpha}, {\boldsymbol \phi}, {\boldsymbol \pi_{Q}}) $,
is dependent  not solely on  $ {\boldsymbol \pi_{Q}}$,
but also on $  {\boldsymbol \phi}$ and ${\boldsymbol\alpha} $.
In the preceding subsection,
we have devised a DNN estimator,  $\widehat {\boldsymbol \phi} $,
  for ${\boldsymbol \phi} $.
A straightforward estimator  for ${\boldsymbol\alpha}$
is derived as follows:
$\widehat{\boldsymbol\alpha} = (\log ( \hat{\pi}_{P,k}/\hat{\pi}_{P,1}), \cdots,
 \log( \hat{\pi}_{P,k}/\hat{\pi}_{P,k-1} ))^{\top}$,
 where  $\hat{\pi}_{P,l} = (1/n_{P})\sum_{j=1}^{n_{P}} I(Y_{j}^{P}=l)$ for $l \in [k]$.
Substituting  $ (\widehat{\boldsymbol\alpha},  \hat {\boldsymbol \phi}) $
for  $ ( {\boldsymbol\alpha}, {\boldsymbol \phi}) $,
we obtain a pseudo-likelihood function of $\bPiQ$, i.e.
$\sum_{j=1}^{n_{Q}}\ell_{Q}(X_{j}^{Q};\widehat{\boldsymbol\alpha}, \widehat{\boldsymbol \phi}, {\boldsymbol \pi_{Q}})$.
To estimate   $\bPiQ $, we propose utilizing
  the pseudo maximum likelihood estimator (PMLE),
\ba
\label{definie-PMLE-multi}
\widehat{\boldsymbol \pi}_{Q}  =
\operatorname*{argmax}_{\bPiQ\in\Pi}\sum\limits_{j=1}^{n_{Q}}
\ell_{Q}(X_{j}^{Q};\widehat{\boldsymbol\alpha}, \widehat{\boldsymbol \phi}, {\boldsymbol \pi_{Q}}),
\ea
where   $\Pi:=\{\bPiQ | \bPiQ  \in[0,1]^{k-1}, 1-\bPiQ^{\top}\mathbf{1}_{k-1}\geq 0\}$,
and
$\mathbf{1}_{k-1}$ denotes a $(k-1)$-dimensional vector of $1$.
This estimation procedure for ${\boldsymbol \pi}_{Q}$ is
exceptionally versatile and adaptable
as it works once reasonable estimates of the density ratios $\{f_{l}(\cdot)/f_{k}(\cdot)\}_{l=1}^{k-1}$ or equivalently $\{P_{X\mid Y}(x\mid l)/P_{X\mid Y}(x\mid k)\}_{l=1}^{k-1}$ are available.

\subsection{The proposed classifier}

Ideally the optimal classifier for the target domain
in terms of excess-risk is the Bayes classifier, defined as
$f_Q^*(x) = \arg\max_{y\in [k]}  \eta_{Q, y}(x)$.
To estimate this classifier, we propose  estimating
the regression function $\eta_{Q, y}(x)$.
For $y \in [k]$,   Bayes's formula dictates that
$
\eta_{Q,y}(x)= \pi_{Q,y}Q_{X\mid Y}(x\mid y) \big/ \sum_{j=1}^{k}\pi_{Q,j}Q_{X\mid Y}(x\mid j).
$
By leveraging the GCS model \eqref{GCS}, we can further derive
\ba
\label{equation-regression-target}
\eta_{Q,l}(x) = \frac{\pi_{Q,l}P_{X\mid Y}(x\mid l)\cdot h(x)}{\sum\limits_{j=1}^{k}\pi_{Q,j}P_{X\mid Y}(x\mid j)\cdot h(x)}=\frac{\eta_{P,l}(x)  \pi_{Q,l}/\pi_{P,l} }
{\sum\limits_{j=1}^{k} \eta_{P,j}(x) \pi_{Q,j}/\pi_{P,j}},
\ea
where the last equality is a consequence of  \eqref{density-ratio-pmle}.

By substituting $\hat\eta_{P,j}(x), \hat \pi_{Q,j}$ and $\hat\pi_{P,j}$ into   equation \eqref{equation-regression-target},
we have the following   plug-in estimators
\ba
\label{defin-classfier-proposed}
\hat{\eta}_{Q,l}(x) = \frac{\hat{\eta}_{P,l}(x)
 \hat{\pi}_{Q,l}/\hat{\pi}_{P,l}
}{\sum\limits_{j=1}^{k}\hat{\eta}_{P,j}(x) \hat{\pi}_{Q,j}/\hat{\pi}_{P,j}  },
\quad l=1,2,\cdots,k-1.
\ea
Ultimately,  our proposed   classifier for  unlabelled target data  is given by
\ba
\label{defin-classifier-multi}
\hat{f}(x) = \operatorname*{argmax}_{l\in [k]}\hat{\eta}_{Q, l}(x)= \operatorname*{argmax}_{l\in [k]}\frac{\hat{\pi}_{Q, l}}{\hat{\pi}_{P,l}}\hat{\eta}_{P, l}(x).
\ea
A notable advantage of our classifier is that it circumvents the problem of estimating
 the function $h(\cdot)$  in the GCS model \eqref{GCS}.
This is because  $h(\cdot)$ appears in both the numerator and the denominator of equality \eqref{equation-regression-target},
 resulting in its  cancellation.
Consequently,  it is unnecessary to estimate $h(\cdot)$ in  \eqref{defin-classfier-proposed}
although $h(\cdot)$ is theoretically identifiable
(as stated in Lemma \ref{prop-identifiability-GCS}) and can be consistently estimated.

\section{Asymptotics}
\label{sec-theory}

This section investigates the asymptotic properties of
the proposed estimators $\hat {\bm \eta}_P$ and $\hatbPiQ$,
as well as the proposed classifier   $\hat{f}(x)$ when utilized with target data.

\subsection{Concentration bounds for $\hat {\bm \eta}_P$ and $\hatbPiQ$}

We commence by  presenting concentration bounds for the DNN estimator $\hat {\bm \eta}_P$
of ${\bm \eta}_P$ and  the PMLE $\hat{\bm \pi}_Q $ of $\bPiQ $.
To facilitate the ensuing discussion,
we introduce additional  notations.
For any two sequences $(a_{n_{P}})_{n_{P}}$ and $(b_{n_{P}})_{n_{P}}$,
we denote $a_{n_{P}}\lesssim  b_{n_{P}}$ if there exists a positive constant
$C$ such that $a_{n_{P}}\leq Cb_{n_{P}}$ for all $n_{P}$.
The notation $a_{n_{P}}\asymp b_{n_{P}}$
signifies that both  $a_{n_{P}}\lesssim  b_{n_{P}}$ and $b_{n_{P}}\lesssim  a_{n_{P}}$.
We define $\bxi^{0}(x)  =
 (\alpha_{1}^{0}+\phi_{1}^{0}(x), \cdots, \alpha_{k-1}^{0}+\phi_{k-1}^{0}(x) )^{\top}$
 and  introduce  $\boldsymbol{\chi}  =  \mathbb{E}_{Q_{X}}
 [(\exp\{\bxi^{0}(X)\}-\mathbf{1}_{k-1})^{\otimes 2}/ (\mathbf{1}_{k-1}^{\top}\exp\{\bxi^{0}(X)\}+1 )^{2} ]$,
 where $A^{\otimes2} = AA^\T$ for any vector or matrix $A$.
 Furthermore,  we denote $c_{\lambda} =  \lambda_{\text{min}}(\boldsymbol{\chi})$.
Throughout our discussion,
  $\mathbb{E}_{\mathcal{D}}$ will be used
  to represent the expectation taken with respect to the randomness of $\mathcal{D}$.

Drawing inspiration from \cite{Schmidt-Hieber2020} and \cite{Zhong2022},  we impose
 the following assumption  pertaining to the network architectures  of $\mathcal{F}(K,s,{\bf p}, D)$.
\begin{assumption}
	\label{assum-network-structure} The neural network $\mathcal{F}(K,s,{\bf p}, D)$ satisfies
(i) $D\geq  (M \vee 1)$,
	(ii) $\sum\limits_{i=0}^{q}\log_{2}(4\tilde{k}_{i}\vee4\theta_{i})\log_{2}n_{P}\leq K\lesssim \log n_{P}$,
	(iii) $n_{P}\gamma_{n_{P}}^{2}\lesssim  \min_{i \in [K]}p_{i}\leq \max_{i \in [K]}p_{i}\leq n_{P}$,
	(iv) $s \asymp n_{P} \gamma_{n_{P}}^{2} \log n_{P}$.
\end{assumption}

For any function $\boldsymbol{f}_{1}$ and $\boldsymbol{f}_{2}$ from
$[0,1]^{d}$ to $\mathbb{R}^{k-1}$, let $d_{P}^{2}(\boldsymbol{f}_{1},
\boldsymbol{f}_{2}) = \mathbb{E}_{P_{X}}\{ \|\boldsymbol{f}_{1}(X)-\boldsymbol{f}_{2}(X)\|^{2}\}$.
Theorem \ref{thm-rate-dnn}   establishes  a concentration bound for the DNN estimator $\hatbetap$.

\begin{theorem}
	\label{thm-rate-dnn}
	Suppose the data $\mathcal{D}$ is generated according to the distribution pair $\mathcal{P}$,
and  that  Assumption \ref{assume-nonparametric-setting} (ii)  and Assumption \ref{assum-network-structure}
are satisfied.
Then, for any $\varepsilon\in(0,1)$ and sufficiently large $n_{P}$,
it holds with probability at least $1-\varepsilon$ that
	\bas
	d_{P}\left(\hatbetap, \betapt\right)\leq \zeta_{1}\sqrt{\log \left(\frac{\zeta_{2}}{\varepsilon}\right)} \cdot \gamma_{n_{P}}\log ^{2}n_{P},
	\eas
	where $\zeta_{1}$ and $\zeta_{2}$ are
positive constants depending solely on $M$, $D$ and $k$.
\end{theorem}

According to Theorem \ref{thm-rate-dnn},
the convergence rate of our DNN estimator of ${\bm\eta}_P$
is of the order $  \gamma_{n_{P}}\log^{2}n_{P} $, which hinges on the intrinsic dimension  $\tilde{k}_{i^{\star}}$. The convergence rate, $\gamma_{n_{P}}\log ^{2}n_{P}$,
 aligns with that of the estimator  presented in \cite{Schmidt-Hieber2020}
 in the realm of nonparametric regression with
 Gaussian errors,
 under the composite smoothness assumption  on the true regression function. When $k=2$, the classification function $\eta_{P, 1}$  can be regarded as
the propensity score function in causal inference and missing data problems.
It serves as a  crucial component in the widely-used  augmented inverse
probability weighting  \citep[AIPW]{robins1994estimation}  and
double machine learning  methods \citep[DML]{Chernozhukov2018}.
Along with other conditions,  the convergence rate result  of our DNN estimator $\hat \eta_{P, 1}$
presented in Theorem \ref{thm-rate-dnn} can be used to verify
whether  the AIPW and DNN methods achieve the  asymptotically semiparametric efficiency lower bound \citep{Chernozhukov2018}. 
%{\blue Theorem \ref{thm-rate-dnn} holds promise for performing statistical inference on the treatment effects in casual inference and missing data, typically necessitating a model for the propensity score, for example, the inverse propensity score matching \citep{hirano2003efficient} or the double robust methods \citep{robins1994estimation}.}

To derive a  concentration bound for $\hatbPiQ$,
we necessitate  an additional
restriction  on the neural networks within $\mathcal{F}(K,s,{\bf p}, D)$.
Specifically,  the constant $s$ governs
the   number of active parameters or  the sparsity of these networks.

\begin{assumption}
	\label{assumption-piq}  $s = o(n_{Q}\gamma_{n_{P}}+\sqrt{n_{Q}}\log^{-3}n_{P})$ as $n_{P}$ and $n_{Q}$ go to infinity.
\end{assumption}

 Assumption \ref{assumption-piq}  stipulates  that
the  sparsity of  the neural network should  not increase too rapidly.
Provided that this assumption is met,
Theorem \ref{thm-rate-piq} furnishes   a concentration bound for
the proposed PMLE  $\hatbPiQ$, which represents a novel contribution
to the existing literature.

\begin{theorem}
		\label{thm-rate-piq}
		Suppose that the data $\mathcal{D}$ is generated
 according to the distribution  pair $\mathcal{P}$,
 and  that Assumptions \ref{assume-nonparametric-setting}--\ref{assumption-piq} are satisfied.
 Then, for any $\varepsilon\in(0,1)$ and sufficiently large $n_{P}$ and $n_{Q}$,
 it holds   with probability at least $1-\varepsilon$ that
	\ba
	\label{bound-piq}
 \|\hatbPiQ-\bPiQ^{0} \|\leq \zeta_{3}\sqrt{\log\left(\frac{\zeta_{4}}{\varepsilon}\right)} \cdot (\gamma_{n_{P}}\log^{2}n_{P}+n_{Q}^{-1/2} ),
	\ea
	where $\zeta_{3}$ and $\zeta_{4}$ are some positive constants depending solely on $\tilde{c}$, $M$, $D$, $k$ and $c_{\lambda}$.
\end{theorem}

Theorem \ref{thm-rate-piq} suggests, in broad sense, that
the convergence rate of our class proportion estimator
is of the order $  (\gamma_{n_{P}}\log^{2}n_{P}+n_{Q}^{-1/2} )$.
When the size of the target data $n_Q$ is substantial,
the term  $\gamma_{n_{P}}\log^{2}n_{P}$ becomes the dominant factor,
whereas for smaller $n_Q$,  the term $n_{Q}^{-1/2}$
prevails.
This finding presents a more intuitive and streamline form
compared to the result in    \cite{iyer2014maximum} [e.g.,Theorem 1].
%{\blue The results in Theorem \ref{thm-rate-piq} immediately imply that the convergence rates for the estimation of the shift function $\pi_{Q,l}/\pi_{P,l}$ mirror those specified  in Theorem \ref{thm-rate-piq}, providing a robust theoretical basis for deriving the generalization bound for the minimizer of the weighted empirical risk within the context of label shift, as illustrated in \citep{ azizzadenesheli2019regularized}}.
Theorem \ref{thm-rate-piq} indicates that the convergence rate of $\hat \pi_{Q, l}$ is no faster than $n_P^{-1/2}$. Because the convergence rate of $\hat \pi_{P, l}$  is $n_P^{-1/2}$, it also indicates that the estimate  $\hat \pi_{Q,l}/\hat \pi_{P,l}$ has the same convergence rate as $\hat \pi_{Q, l}$. This convergence rate establishes a theoretical cornerstone for deriving the generalization error bound of the minimizer of the weighted empirical risk under label shift \citep{ azizzadenesheli2019regularized}. More importantly, the results  presented in  Theorems \ref{thm-rate-dnn}--\ref{thm-rate-piq}
establish the groundwork for achieving
the minimax optimal convergence rate of
our proposed classifier in terms of excess-risk.

\subsection{Minimax optimality of the proposed classifier}

We initially derive an upper bound for  the excess-risk $\mathcal{E}_Q (\hat{f} )$
of the proposed classifier  $\hat{f}(\cdot)$, as  defined in \eqref{defin-classifier-multi}.

\begin{theorem}
\label{thm-rate-excess-risk}
Suppose that the data $\mathcal{D}$ is generated according to the distribution pair $\mathcal{P}$,
 and that   Assumptions \ref{assume-nonparametric-setting}--\ref{assumption-piq} are satisfied.
 Then, for any $\varepsilon\in(0,1)$ and sufficiently large $n_{P}$ and $n_{Q}$,
 it holds with   probability at least $1-\varepsilon$ that
\ba
\label{bound-excess-risk}
\mathcal{E}_Q (\hat{f} )\leq \zeta_{5}\sqrt{\log\left(\frac{\zeta_{6}}{\varepsilon}\right)}
\cdot \left(\gamma_{n_{P}}\log^{2}n_{P}+n_{Q}^{-1/2}\right),
\ea
	where $\zeta_{5}$ and $\zeta_{6}$ are some positive constants depending solely on $\tilde{c}$, $M$, $D$, $k$ and $c_{\lambda}$.	
\end{theorem}

Given that for  a nonnegative random variable $T$,
$\e T  = \int_0^{\infty} P(T\geq t) dt$,
  Theorem \ref{thm-rate-excess-risk} implies
\ba
\label{f-hat-upper-bound}
\e_{\mathcal{D}}\{ \mathcal{E}_Q(\hat{f})\}
\leq \sqrt{2\pi}\cdot \zeta_6  \zeta_5  \left(\gamma_{n_{P}}\log^{2}n_{P}+n_{Q}^{-1/2}\right),
\ea
indicating that  the expected excess-risk  $\e_{\mathcal{D}}\{ \mathcal{E}_Q(\hat{f})\}$
shares the same  asymptotic convergence rate $(\gamma_{n_{P}}\log^{2}n_{P}+n_{Q}^{-1/2})$ as
the  excess-risk $\mathcal{E}_Q (\hat{f} )$.

Furthermore, Theorem \ref{thm-rate-excess-risk} elucidates   that
the contributions of the source and target data to
the convergence rate of the proposed classifier
are $\gamma_{n_{P}}\log^{2}n_{P}$ and $n_{Q}^{-1/2}$, respectively.
Notably, the converge rate of our classifier is independent of $d$,
allowing the proposed method to alleviate
the curse of dimensionality when
the dimension $d$ of $X$ is high but  the intrinsic dimension
of the composite function $\bphit(\cdot)$  remains relatively low.
To illustrate this,
we present two examples showcasing the selection of $\bphit(\cdot)$
with low-dimensional structures.

\begin{example}
	[Generalized additive model, GAM  \citep{generalizedadditive2007}]
	\label{example-GAM}
Consider functions
  $\phi_{1,j}:[0,1]\mapsto[-M_{0}, M_{0}]$ for some $M_{0}>0$ and  $1\leq j\leq d$,
and  $\tilde{H}(\cdot):\mathbb{R}\mapsto\mathbb{R}$.
Define   $\phi_{1}^{0}(x) = \tilde{H}(\sum_{j=1}^{d}\phi_{1,j}(x_{(j)}))$,
where   $x_{(j)}$  denotes the $j$-th component of $x$ and $d$ is a positive integer.
If we take $g_{0}(x) = (\phi_{1,1}(x_{(1)}),\cdots, \phi_{1,d}(x_{(d)}))^{\top}$,
$g_{1}(x) = \sum_{j=1}^{d}x_{(j)}$ and $g_{2}(\cdot) = \tilde{H}(\cdot)$, then  $\phi_{1}^{0}(x)= g_{2}\circ g_{1}\circ g_{0}(x)$.
 Assume that $\phi_{1,j}\in\mathcal{H}_{1}^{\theta_{0}}([0,1], M_{0})$ for $j \in [d]$ and some $\theta_{0}>0$, and $\tilde{H}\in\mathcal{H}_{1}^{\theta_{2}}([-dM_{0}, dM_{0}], M_1)$ for some $\theta_{2}>0$ and $M_1>0$.
 It is clear that   $g_{1}(\cdot)\in\mathcal{H}_{d}^{\theta_{1}}([-M_{0}, M_{0}]^{d}, dM_{0}+1)$  for any $\theta_{1}>1$.
Therefore
$
	\phi_{1}^{0}\in\mathcal{G}(2, (\theta_{0}, (\theta_{0}\vee 2)d,\theta_{2}), (d, d, 1,1), (1, d,1), (dM_{0}+1)\vee M_1).
$
According to Theorem \ref{thm-rate-excess-risk},
the contribution of the source data
to the  convergence rate of   $\mathcal{E}_Q (\hat{f} )$
  is of the order
$
 \big( n_{P}^{-\frac{\theta_{0}(\theta_{2}\wedge 1)}{2\theta_{0}(\theta_{2}\wedge 1)+1}}+n_{P}^{-\frac{\theta_{2}}{2\theta_{2}+1}}
 \big)\log^{2}(n_{P}).
$
Here are some important observations:

 \bit
 \item
 When $\theta_{0}=\theta_{2}\geq 1$, the rate simplifies to $n_{P}^{-\frac{\theta_{2}}{2\theta_{2}+1}}\log^{2}(n_{P})$,
 which, up to a logarithm factor,
 aligns with the rate in Theorem 2.1 of \cite{generalizedadditive2007}
 in the context of penalized   nonparametric (generalized additive) regression.

 \item	
The GAM reduces to an additive model  \citep{stone1985additive} when $\tilde{H}$ is
the identity function.
 In this situation,  the convergence rate contribution of the source data is
 of the order $n_{P}^{-\frac{\theta_{0}}{2\theta_{0}+1}}\log^{2}(n_{P})$,
 which, up to a logarithm factor,  coincides with the optimal minimax rate in Corollary 1
 of \cite{stone1985additive} for additive (mean) regression.

 \item
 When $\phi_{1,j}$ is a linear function of $x_{(j)}$ for $j\in [d]$,
 the GAM reduces to the single index model \citep{Ichimura1993semiparametric},
 and the convergence rate contribution of the source data is of the order $n_{P}^{-\frac{\theta_{2}}{2\theta_{2}+1}}\log^{2}(n_{P})$,
 which, up to a logarithm factor,  matches the optimal rate in Theorem 3 of \cite{Ma2016inference}
 for single-index (quantile) regression.
\eit
In  all three cases, the convergence rate is independent of   the
  dimension $d $ of $X$, which can be any positive integer.
This underscores that the ability of our method  to  mitigate the  curse-of-dimensionality.

\end{example}

\begin{example}[Generalized hierarchical
	interaction model  \citep{kohler2016nonparametric}]
	\label{example-GHIM}
	For a fixed positive integer $d^{\star}\in \mathbb{N}_{+}$ ($d^{\star}<d$),
consider vectors $\beta_{i}^{1},\cdots,\beta_{i}^{d^{\star}}\in\mathbb{R}^{d}$
($i\in [d^{\star}]$),
functions $\tilde{f}_{1},\cdots, \tilde{f}_{ d^{\star}}:\mathbb{R}^{d^{\star}}\mapsto\mathbb{R}$  and a function $\tilde{H}:\mathbb{R}^{d^{\star}}\mapsto\mathbb{R}$  such that
	\bas
\phi_{1}^{0}(x) =  \tilde{H}\left(\tilde{f}_{1}(x),\cdots \tilde{f}_{ d^{\star}}(x)\right),
 ~\tilde{f}_{i}(x) = \tilde{f}_{i}(x^{\top}\beta_{i}^{1},\cdots, x^{\top}\beta_{i}^{d^{\star}}), ~~ 1\leq i\leq d^{\star}.
	\eas
Define the functions:
$ g_{0}(x) = (x^{\top}\beta_{1}^{1},\cdots, x^{\top}\beta_{1}^{d^{\star}}, \cdots, x^{\top}\beta_{d^{\star}}^{1},\cdots, x^{\top}\beta_{d^{\star}}^{d^{\star}})^{\top}$,
$
	g_{1}(t_{1}, \cdots,t_{d^{\star}})
=  (\tilde{f}_{1}(t_{1}),\cdots,\tilde{f}_{d^{\star}}(t_{d^{\star}}))^{\top},
$
	and $g_{2}(t) = \tilde{H}(t)$, where   $t, t_{1},\cdots, t_{d^{\star}}$ are  $d^{\star}$-dimensional  vectors.
Assume that $\|\beta_{i}^{j}\|\leq M_{0}$ for all $ i, j\in [d^{\star}]$ for some $M_{0}>0$.
Consequently, for any $\theta_{0}>1$, the linear functions
 $x^{\top}\beta_{i}^{j}$ belong to $\mathcal{H}_{d}^{\theta_{0}}([0,1]^{d}, dM_{0})$.
 Furthermore, suppose that   $\tilde{f}_{1},\cdots, \tilde{f}_{d^{\star}}, \tilde{H}$ all belong to $\mathcal{H}_{d^{\star}}^{\theta}([-M_{1}, M_{1}]^{d^{\star}}, M_{1})$ for some $\theta>0$ and sufficiently large $M_{1}$.
Then it follows that
	\bas
	\phi_{1}^{0}\in\mathcal{G}\left(2, ((\theta\vee 2)d/(\theta\wedge 1)^{2}, \theta, \theta), (d, (d^{\star})^{2}, d^{\star}, 1), (d, d^{\star}, d^{\star}), M_{1}\vee dM_{0}\right).
	\eas
According to Theorem \ref{thm-rate-excess-risk},
the  contribution of the source data to the convergence rate  of
$\mathcal{E}_Q (\hat{f} )$ is of the order
$ \big(n_{P}^{-\frac{\theta(\theta\wedge 1)}{2\theta(\theta\wedge 1)+d^{\star}}}+n_{P}^{-\frac{\theta}{2\theta+d^{\star}}}\big)\log^{2}(n_{P}).
$
 When $\theta\geq 1$, the rate simplifies to $n_{P}^{-\frac{\theta}{2\theta+d^{\star}}}\log^{2}(n_{P})$,
 which aligns with the rate in Theorem 1 of \cite{BauerandKohler2019} in the context of nonparametric regression, up to a logarithm factor.
Notably,  the curse of dimensionality  is  mitigated
 here because  $ d^{\star}$ is strictly smaller than $d$.

\end{example}

Next, we establish a lower bound for the   excess-risks of all classifiers from $\mathcal{S}$ to $\mathcal{H}$.
Recall that $i^{\star}\in\operatorname*{argmin}_{i\in \{0,\cdots, q\}}\tilde{\theta}_{i}/(2\tilde{\theta}_{i}+\tilde{k}_{i})$.

\begin{theorem}
	\label{thm-minimax}
	Suppose that the data $\mathcal{D}$ is generated according to the distribution pair $\mathcal{P}$,
and that the source data and target data in $\mathcal{D}$ are independent of each other.
Under the distribution class $\mathcal{P}^{\star}$ such that
$\tilde{k}_{i^{\star}}\leq \min\{k_{0},\cdots, k_{i^{\star}-1}\}$,  there exits a positive constant $\zeta_{7}$  depending solely on $(\tilde{\theta}_{i^{\star}}, \tilde{k}_{i^{\star}},k)$ such that
\ba
	\label{lower-bound-excess-risk}
\inf _{\mathcal{A}:\mathcal{S}\mapsto\mathcal{H}}\sup _{\mathcal{P} \in \mathcal{P}^{\star}} \mathbb{E}_{\mathcal{D}}\left\{
\mathcal{E}_Q\left(\mathcal{A}\left(\mathcal{D}\right)\right)\right\}
\geq\zeta_{7}\left( \gamma_{n_{P}}+n_{Q}^{-1/2} \right).
\ea
	
\end{theorem}

The lower bound of the  expected excess-risk,
as established in Theorem \ref{thm-minimax},
is contingent upon the condition
$\tilde{k}_{i^{\star}}\leq \min\{k_{0},\cdots, k_{i^{\star}-1}\}$. This implies that no additional
dimensions are incorporated at the intrinsic dimension level
within the composite function, thereby precluding
the scenario where $\tilde{k}_{i^{\star}}$ surpasses
the covariate dimension $d$.
Notably, the lower bound result in Theorem \ref{thm-minimax}
retains its validity when the condition is fortified  to
 $\tilde{k}_{i}\leq \min\{k_{0},\cdots, k_{i-1}\}$ for all $1\leq i \leq k$  (See, e.g., \cite{Schmidt-Hieber2020}).  This signifies that
no dimensions are added on deeper abstraction levels
in the composition of functions.

According to Theorem  \ref{thm-minimax},
the convergence rates of the expected excess-risk for all classifiers
transitioning
 from $\mathcal{S}$ to $\mathcal{H}$
 do not surpass   $\gamma_{n_{P}} +n_{Q}^{-1/2}$.
In the meantime,    Theorem  \ref{thm-rate-excess-risk}
  indicates that
  the convergence rate of the expected excess-risk for the proposed
classifier $\hat f$ is no slower than $\gamma_{n_{P}}\log^{2}n_{P}+n_{Q}^{-1/2}$;
see inequality \eqref{f-hat-upper-bound}.
Combining these results, we immediately deduce the   minimax optimal rate
for all classifiers
and establish the minimax optimality of the proposed classifier $\hat f$.

 \begin{coro}
	\label{thm-minimax-f-hat}
Assume the conditions in Theorems  \ref{thm-rate-excess-risk} and \ref{thm-minimax}.
Then, the minimax optimal rate of convergence for all classifiers
 from $\mathcal{S}$ to $\mathcal{H}$ is  $\gamma_{n_{P}} +n_{Q}^{-1/2}$
 up to a logarithm factor $\log^2 (n_P)$,
and  the proposed classifier $\hat f$ attains this optimal rate  asymptotically.
	
\end{coro}

\begin{remark} \label{remark1}
In nonparametric binary classification problems, a commonly employed   condition is the marginal assumption, which posits the existence of
 constants $C_0>0$ and $\alpha\geq 0$ such that
 either $Q( 0<|\eta_{Q,1}^{0}(X) - 0.5| \leq t) \leq C_0t^{\alpha}$ \citep{Fastlearning2007,maity2022minimax}
or $Q(  |\eta_{Q,1}^{0}(X) - 0.5| < t) \leq C_0t^{\alpha}$ \citep{TonyCai2021}
holds for all $t>0$.
However, in our nonparametric classification setting,
as outlined in Assumption \eqref{assume-nonparametric-setting},
we refrain from employing such an assumption for two primary reasons.
Firstly, the margin assumption can often be violated.
For example, in   binary classification scenarios,
this assumption fails when
  $\eta_{Q,1}^{0}(X)$ equals 0.5 with a positive probability.
Secondly,   verifying the margin  assumption
is challenging due to the unavailability of the target data labels.
When  the optimal rates derived
in \cite{maity2022minimax, maity2024linear} and \cite{TonyCai2021}
become more favorable as the parameter  $\alpha$ in the margin assumption
increases (indicating a stronger margin assumption),
the rates still depend on the covariate dimension $d$.
Consequently,   their methods  remain susceptible to
the curse of dimensionality.
\end{remark}

In the context of binary classification, \cite{maity2022minimax}   derived a lower bound for the expected excess-risk of all classifiers  using  labelled source data and unlabelled target data. Specifically, when the margin assumption is violated (i.e., $\alpha=0$),
 the optimal rate obtained in  \cite{maity2022minimax} [Theorem 12] is  of the order
$n_{P}^{\frac{-\tilde{\beta}}{2\tilde{\beta}+d}}+n_{Q}^{-1/2}$,
where $\tilde{\beta}$ represents the   smoothness index of the probability density functions $f_{1}(X)$ and $f_{2}(X)$.
Notably, in scenarios where the intrinsic dimension of the composite function $\phi_{1}^{0}$ is relatively low (as illustrated in   Examples \ref{example-GAM}--\ref{example-GHIM}), the   convergence rates presented in Theorems \ref{thm-rate-excess-risk}--\ref{thm-minimax}
surpass those of  \cite{maity2022minimax} when $d$ is   large.
When $\phi_{1}^{0}$ is fully non-parametric,
 the convergences rate of our classifier
aligns with  that of \cite{maity2022minimax}.
However, under the margin assumption  with $\alpha>0$
 and when $n_{P}=o( n_{Q})$,
 the optimal rate in \cite{maity2022minimax} is of the order
$n_{P}^{\frac{-\tilde{\beta}(1+\alpha)}{2\tilde{\beta}+d}}$,
which still suffers from  the curse of dimensionality.
In contrast, our classifier's convergence rate    is of the order $\gamma_{n_{P}}$,
which is  faster than that of \cite{maity2022minimax}
when the intrinsic dimension of $\phi_{1}^{0}$ is low and $d$ is large.

\section{Simulations}
\label{sec-simu}
\label{subsection-simu-method}

\subsection{Methods under comparison}

We  conducted  simulation experiments to
assess the finite-sample performance of the proposed
DNN and PMLE-based Classifier (DNN-PC) defined in \eqref{defin-classifier-multi}.
Although our DNN-PC is versatile and can be applied to  multi-class classification problems,
in this simulation study,  we  focused specifically
on the binary-classification setting,  where $k=2$.
The Bayes classifier for the target data is  $ {f}(x) = 2-  I(\eta_{Q,1}(x)\geq 1/2)$.
For comparative purposes,
we  also include several popular competitors, which are
variations of  the Bayes classifier utilizing different estimates of $\eta_{Q,1}(x)$.

\bit
\item
Maity-PC:  The classifier proposed by \cite{maity2022minimax}, employs a PMLE.
Specifically,   $\hat{P}_{X\mid Y=y}$  represents
the kernel-density estimator  of $ {P}_{X\mid Y=y}$  for $y=1, 2$,
and
$\tilde \pi_{Q, 1}$ be  the PMLE of  $\pi_{Q, 1}$  obtained
by replacing the estimates of  $ {P}_{X\mid Y=1}/{P}_{X\mid Y=2}$
with the ratio of their kernel-density estimators.
Maity-PC  estimates $\eta_{Q,1}(x)$ by
\bas
\tilde{\eta}_{Q,1}(x)&=&
\frac{ \tilde \pi_{Q, 1} \hat{P}_{X\mid Y }(x\mid 1)}{ \tilde \pi_{Q, 1}
\hat{P}_{X\mid Y }(x\mid 1)+(1-\tilde \pi_{Q, 1})\hat{P}_{X\mid Y }(x\mid 2)}.
\eas

\item{Saerens}: This classifier, proposed by
\cite{Saerens2002adjusting}, estimates   $\eta_{P,1}$   using
the  $\text{K}$-nearest neighbors (KNN) method
and $\pi_{Q, 1}$ using  an Expectation-Maximization (EM) algorithm.

\item
Maity-IC:
This is the ideal classifier proposed by
\cite{maity2022minimax}, which is essentially
 Maity-PC with  $\tilde \pi_{Q, 1}$ replaced by
 the ideal estimator $\hat{\pi}_{Q,\text{ideal}}
 = \sum_{j=1}^{n_{Q}}(2-Y_{j}^{Q})/n_{Q}$.

\item
Maity-LA:
This classifier, proposed by \cite{maity2024linear},
is  based on a KNN estimator for $\eta_{P,1}$ and  Linear Adjustment model, which holds under GCS model.

\eit
The first three classifiers (DNN-PC, Maity-PC, and Saerens) are designed
specifically for unsupervised transfer learning.
The last two (Maity-IC  and Maity-LA) are intended  for supervised transfer learning,
and are not applicable in our settings.
However, they  are  included here solely for comparison purposes.

To implement the proposed DNN-PC method,
we utilize the {\it Adam} optimizer in the {\it Deep Learning Toolbox} of  Matlab
to solve the optimization problem defined in \eqref{definie-NMPLE-multi}.
Subsequently, we calculate the maximizer of the log-partial likelihood presented in \eqref{definie-PMLE-multi}
 by resolving the associated score equation.
This process is facilitated by  the {\it fzero} solver available  in   Matlab.

The methods under comparison  involves numerous hyperparameters.
Specifically, in the proposed DNN-PC method,
the  hyperparameters include  the number of hidden layer,
the number of neurons per  hidden layer,
 the learning rate,  and mini-batch size utilized by
    the {\it Adam} optimizer \citep{kingma2014adam}.
Notably, we maintain a consistent number of neurons across all hidden layers
 in the pertinent DNN.
The implementations of  Maity-IC  and Maity-PC involve kernel-based estimators,
where we adopt the standard   normal density function as
the kernel function and consider  the bandwidth as a hyperparameter.
For the implementations of Maity-LA and Saerens,
we regard the neighborhood size in {\text K-NN} as a hyperparameter.
All these hyperparameters are  selected through grid search.

\subsection{Simulation settings}

We generate $Y^P$ and $Y^Q$ such that $2-Y^{P}$ follows a Bernoulli distribution
with success probability 0.75,  and
$2-Y^{Q}$ follows a Bernoulli distribution with success probability
$ \pi_{Q,1} $,
where   $\pi_{Q,1}$ is considered to be either  0.25
(representing an imbalanced case)
or 0.5 (representing a balanced case).
The covariate $X = (X_{(1)}, X_{(2)}, X_{(3)}, X_{(4)})^{\top}$
is chosen to be a   $4$-dimensional vector.
Given  $Y=y$, the components of $X$, denoted as $X_{(i)}$,
 are independent and identically distributed
in both the source  and target populations.
Three scenarios are considered for the conditional distribution pairs.  Given $Y=y$,
\bit
\item[I.]  Both  $X_{(i)}^P$  and   $X_{(i)}^Q$  follow
$  I(y=1)\times \text{Uniform}(4) + I(y=2)\times \text{Uniform}(2) $,
where $\text{Uniform}(\tilde{o})$ denotes the discrete uniform distribution on $\{1,\cdots,\tilde{o}\}$
for $\tilde{o}\in\mathbb{N}_{+}$;

\item[II.]
Both  both  $X_{(i)}^P$  and   $X_{(i)}^Q$  follow
$ I(y=1)\times \text{Beta}(6, 2) + I(y=2)\times \text{Beta}(2, 6)$,
where $\text{Beta}(a, b)$ is the Beta distribution
with shape parameters $a>0$ and $b>0$;

\item[III.]
$X_{(i)}^P$ follow
$  I(y=1)\times \text{Exp}(e^{0.5}-1) +I(y=2)\times \text{Norm}(1, 1)$ and
$X_{(i)}^Q$ follow $ I(y=1)\times \text{Exp}(1-e^{-0.5}) + I(y=2)\times \text{Norm}(0, 1)$.
Here, $\text{Exp}(\mu_{1})$ is the exponential distribution with mean $\mu_{1}>0$ and $\text{Norm}(\mu_{2}, \sigma)$ denotes
 the normal distribution with mean $\mu_{2}\in\mathbb{R}$ and standard deviation $\sigma>0$.
\eit

The first two scenarios are specifically crafted
to ensure that both the label shift and GCS assumptions are met.
In these scenarios,
the conditional distributions of $X_{(i)}$ given  $Y$ are discrete and continuous, respectively.
However, in  Scenario III,   the label shift assumption is violated,
while  our GCS  assumption remains valid
with $h(x) = \exp(2+\beta^{\top}x)$, where $\beta = (-1,-1,-1,-1)^{\top}$.
Let $p(t|y)$ represent the conditional probability density or mass function
of $X_{(i)}^P$ given $Y^P = y$.
According to the aforementioned data generation settings,
 we have     $\eta_{P,1}( x) = \exp(\phi_{1}(x))/(1+\exp(\phi_{1}(x)))$, where $\phi_{1}(x)  = g_{2}\circ g_{1}\circ g_{0}(x)$ and $g_{2}(t) = \log\{ t/(1-t) \}$, $g_{1}(a,b) =
3a/(3a+b)$ and $g_{0}(x) =  (\prod_{i=1}^{4}p(x_{(i)}\mid 1), \prod_{i=1}^{4}p(x_{(i)}\mid 2) )$.

When generating data,  we set  the sample size pairs $(n_P, n_Q)$ to
$(100\times J, 400)$   with $J\in\{1,2,3,5,6\}$.
To approximate the excess-risk of a generic classifier,
we utilize an independent test sample of size 2500
to compute its empirical excess-risk.
The entire procedure is repeated 200 times, and
the final reported excess-risk of the classifier
is determined by averaging the 200 empirical excess-risks obtained.

\subsection{Simulation results}

\begin{figure}[htbp]
	\centering
	\includegraphics[width=0.9\textwidth]{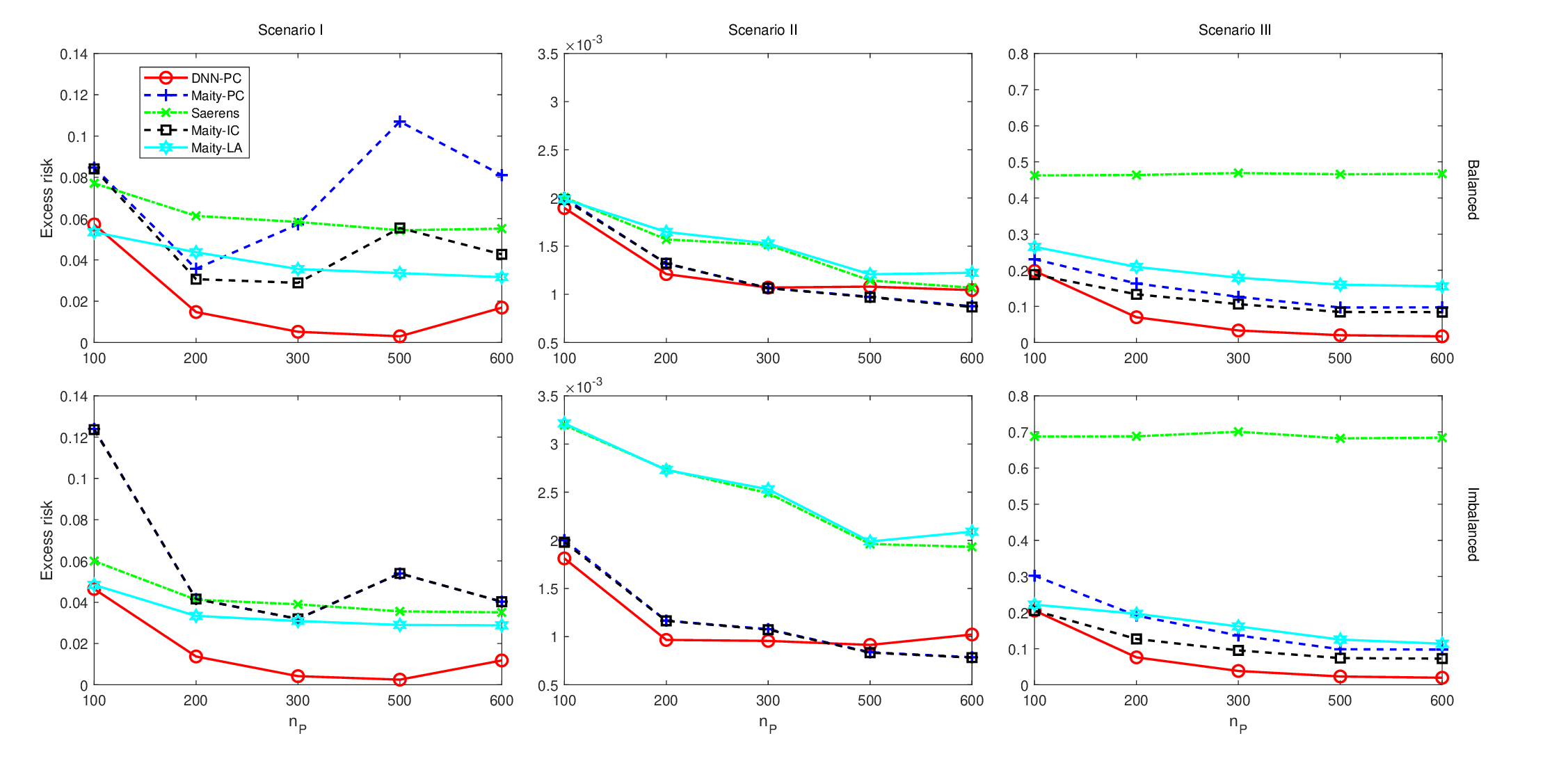}
	\label{fig-excess-risk}
	\caption{
Plots of $n_P$ versus excess-risks of the classifiers under comparison when $n_Q = 400$.
Upper panel:  the  balanced case ($\pi_{Q,1} = 0.5$);
Lower panel:  the  unbalanced case ($\pi_{Q,1} = 0.25$)
}
\end{figure}

Figure  \ref{fig-excess-risk} presents  the   excess-risk results of the  compared classifiers
for different values of $n_P$  and $n_Q = 400$.
Scenarios I and II are crafted to meet the label shift assumption.
In Scenario I,  the covariates have discrete values  and the $\phi_1(\cdot)$ function
is not a  composite smoothness function, thus violating our Assumption \ref{assume-nonparametric-setting}.
Nevertheless,  the proposed DNN-PC classifier consistently exhibits less
 excess-risk compared to all  four competing classifiers, even though Maity-IC and Maity-LA
make the  label information of the target data. A possible reason for this observation is that the kernel and KNN methods underlying  the  four competitors
may have poor performance for discrete covariates.

%  Notably,  Maity-IC and Maity-LA are developed for supervised classification. {\blue For the Maity-PC and Maity-IC methods, their underwhelming performance may stem from the comparatively low efficiency of kernel methods in estimating discrete densities. Regarding the Saerens method, its less than satisfactory performance could be attributed to the absence of convergence guarantees for its EM algorithm. Lastly, despite the Maity-LA method's utilization of labeled target data, the estimation of model parameters within its linear-adjust-model lacks an efficient algorithm for solution.}

In Scenario II, the proposed DNN-PC classifier
 remains  uniformly superior to  Saerens and Maity-LA,
 both of which are based on KNN.
 However, it shows comparable performance to  Maity-PC and Maity-IC, which rely on
the kernel  method.
This is likely due to the covariates being continuous in Scenario II,
for which the kernel density estimator functions effectively.
Additionally,  the competing classifiers necessitate the   estimation  of $P_{X|Y}(x|1), \ldots, P_{X|Y}(x|k)$, while  our DNN-PC classifier only requires the estimation of their ratios, not
the densities themselves. Generally,  estimating the conditional densities $P_{X|Y}(x|i)$ is more difficult than estimating their ratios.
Since the ideal classifier depends solely on the ratios,
by directly estimating the ratios,
our method not only reduces the number of parameters (the baseline conditional density or mass
function) but also yields more reliable results than the competitors.

Scenario III is  designed to violate  the label shift assumption while satisfying
 the GCS assumption. Clearly, in terms of   excess-risk,
the proposed DNN-PC classifier uniformly outperforms the four competitors,
regardless of whether the competitors incorporate
the label information of the target data or not.

\begin{figure}[htbp]
		\label{fig-mse-piQ}
	\centering
		\includegraphics[width=0.9\textwidth]{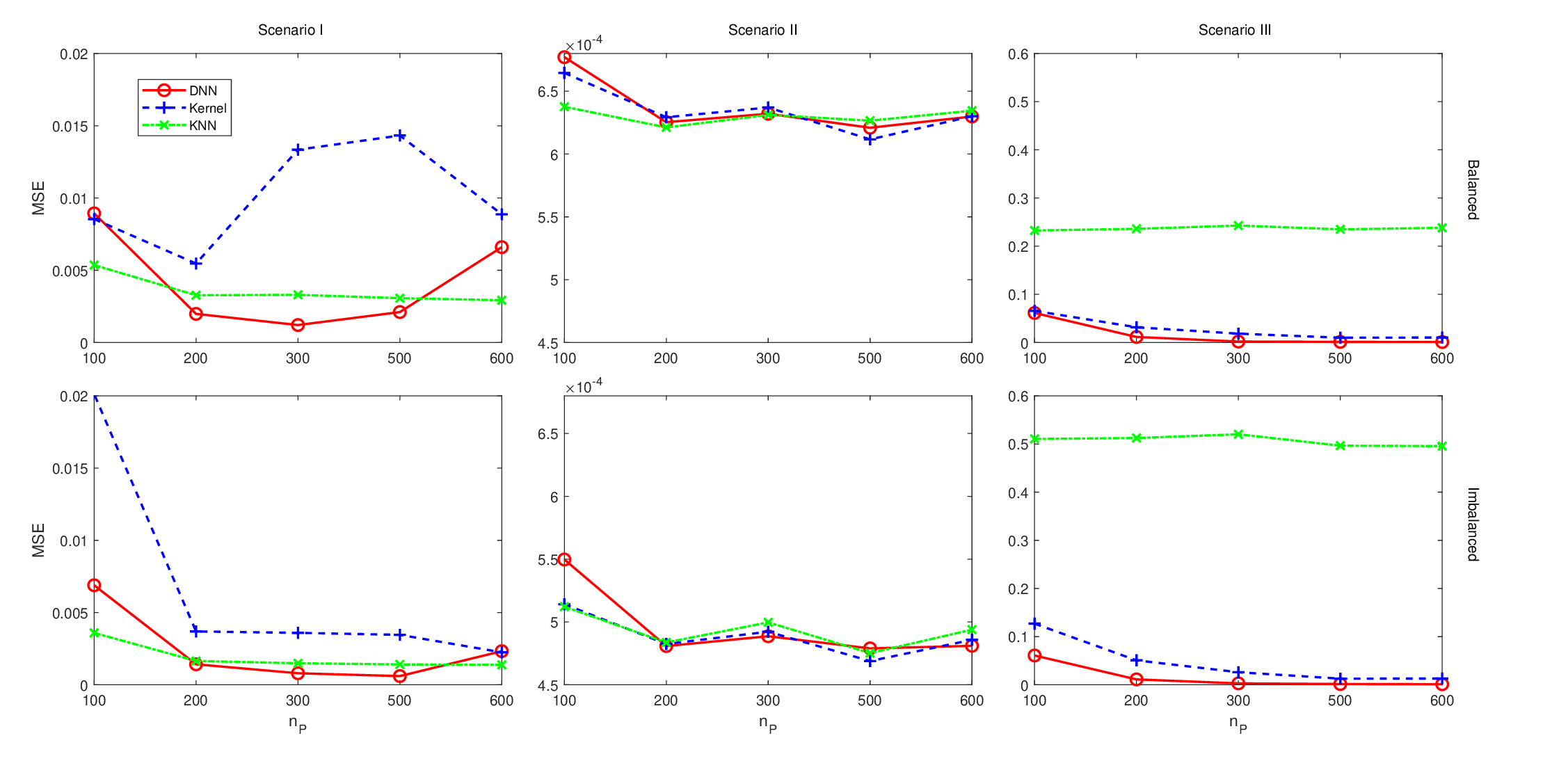}
	\caption{Mean square errors of the three estimators of $\pi_{Q,1}$ under comparison with $n_Q = 400$.
Upper panel:  the  balanced case ($\pi_{Q,1} = 0.5$);
Lower panel:  the  unbalanced case ($\pi_{Q,1} = 0.25$)
}
\end{figure}

The class proportions $\pi_{Q,1}$
are crucial for  the Bayes classifier
in the target population.
Among the classification methods  being compared,
three estimators for $\pi_{Q,1}$ were utilized:
(1) DNN:
The proposed PMLE $\hat{\pi}_{Q,1}$,
 where  the density ratio $P_{X|Y = 1}/P_{X|Y = 2}$
 is estimated  using  the proposed DNN estimators.
(2) Kernel:
The PMLE $\tilde{\pi}_{Q, 1}$
 from Maity-PC,
 with the  density ratio estimate replaced
  by  the ratio of their kernel   estimators,
and (3)
KNN:  The estimator $\hat{\pi}_{Q, \text{EM}}$
from \cite{Saerens2002adjusting}, based on an EM algorithm
after estimating $\eta_{P,1}(\cdot)$ using  K-NN estimators.
Figure \ref{fig-mse-piQ} displays
the mean square error (MSE) results
of these three estimators for $\pi_{Q,1}$
when data were generated according to Scenarios I--III.
In terms of MSE,
the DNN estimator for $\pi_{Q, 1}$ is
comparable to Kernel and KNN in Scenarios I and II,
but consistently outperforms
the latter two in   Scenario III.
Notably, the superiority  of our classifier
often aligns with the superiority of our DNN estimator of $\pi_{Q, 1}$
compared to its competitors.

Additional simulations were conducted with   sample size pairs $(n_P, n_Q)$
set to $  (400, 100\times J )$
for $J\in\{1,2,3,5,6\}$  when data were generated
according to  Scenarios I-III.
The corresponding results are presented in  Section 7 of the supplementary material.
Among the compared methods, our proposed   classifier
 still exhibits the smallest excess-risk in Scenarios I and III.
In Scenario II,  its performance is comparable to the kernel-based classifiers
(Maity-PC and Maity-IC),  and they three uniformly outperform the KNN-based classifiers (Saerens and Maity-LA).
Similarly,  in estimating  $\pi_{Q, 1}$,
our DNN estimator outperforms the kernel and KNN estimators in Scenarios I and III,
and they are comparable in Scenario II.
Unlike the case where  $n_P$ varies,
the excess-risk and MSE remain relatively unchanged as  $n_Q$ increases from 100 to 600
with $n_P$ fixed at 400.

Overall, the proposed estimation method and classifier are comparable or
 superior to   existing competitors
when the label shift assumption is satisfied
and  uniformly outperform them  when
the label shift assumption is violated but GCS is satisfied.
Our method is capable of producing a more precise classifier,
even with limited data availability in the target domain,
while simultaneously demonstrating remarkable generalization ability and efficiency.

\section{Diagnosis of Alzheimer's Disease}
\label{sec-real}

Alzheimer's disease represents a progressive neurological affliction that induces the degeneration and death of brain cells. This process leads to a consistent deterioration in cognitive faculties and functional capabilities.
At present, approximately 50 million individuals globally are living with Alzheimer's disease. Alarmingly, this figure is forecasted to double every five years, with projections indicating that it will surge to 152 million by the year 2050.
The accurate diagnosis of Alzheimer's disease holds paramount significance. It is crucial for effective disease treatment, providing support to patients and their families, and ensuring the proper allocation of medical resources.

In this section, we employ the proposed classification method to analyze
the Alzheimer's Disease Dataset of \cite{rabie_el_kharoua_2024}.
This dataset encapsulates comprehensive health information for 2,149 patients, encompassing demographic specifics, lifestyle factors, medical histories, clinical measurements, cognitive and functional evaluations, symptoms, and Alzheimer's disease diagnoses. We designate the Alzheimer's disease diagnosis for each patient (with 1 indicating a positive diagnosis and 2 indicating a negative one) as the outcome variable $Y$, where $k=2$.  Our analysis includes 18 covariate variables, spanning demographic details (age and gender, comprising 2 variables), lifestyle factors (body mass index, alcohol consumption, physical activity, diet quality, and sleep quality, totaling 5 variables), medical history (family history of Alzheimer's and cardiovascular disease, making up 2 variables), clinical measurements (systolic blood pressure, diastolic blood pressure, total cholesterol levels, low-density lipoprotein cholesterol levels, high-density lipoprotein cholesterol levels, and triglyceride levels, amounting to 6 variables), and cognitive and functional assessments (mini-mental state examination score, functional assessment score, and activities of daily living score, comprising 3 variables).
 The dataset, along with a detailed description of the covariates, is accessible via the URL \url{https://www.kaggle.com/datasets/rabieelkharoua/alzheimers-disease-dataset?resource=download}.

We categorize patients who smoke as belonging to the target domain ($n_{Q}= 620$) and consider non-smoking patients as part of the source domain ($n_{P}= 1529$). Our focus is on training a classification model that can predict whether a patient in the target domain has Alzheimer's disease. Prior to our analysis, all covariates were scaled to fall within the range of 0 to 1.
To assess the validity of the label shift assumption, we introduce an indicator variable $D$, assigning $D=1$ to all target data and $D=0$ to all source data. The label shift assumption translates to the conditional independence of the covariate vector $X$ and the indicator $D,$ given the outcome variable $Y$. We verify this conditional independence using the Invariant Conditional Quantile Prediction method \citep{heinze2018invariant}, available within the R package  {\tt CondIndTests}. The resulting $p$-value of 0.0198 suggests insufficient evidence to uphold the label shift assumption.

\begin{figure}[htbp]
\label{fig-real-test-GCS}
\centering
 \includegraphics[width=14cm, height=6cm]{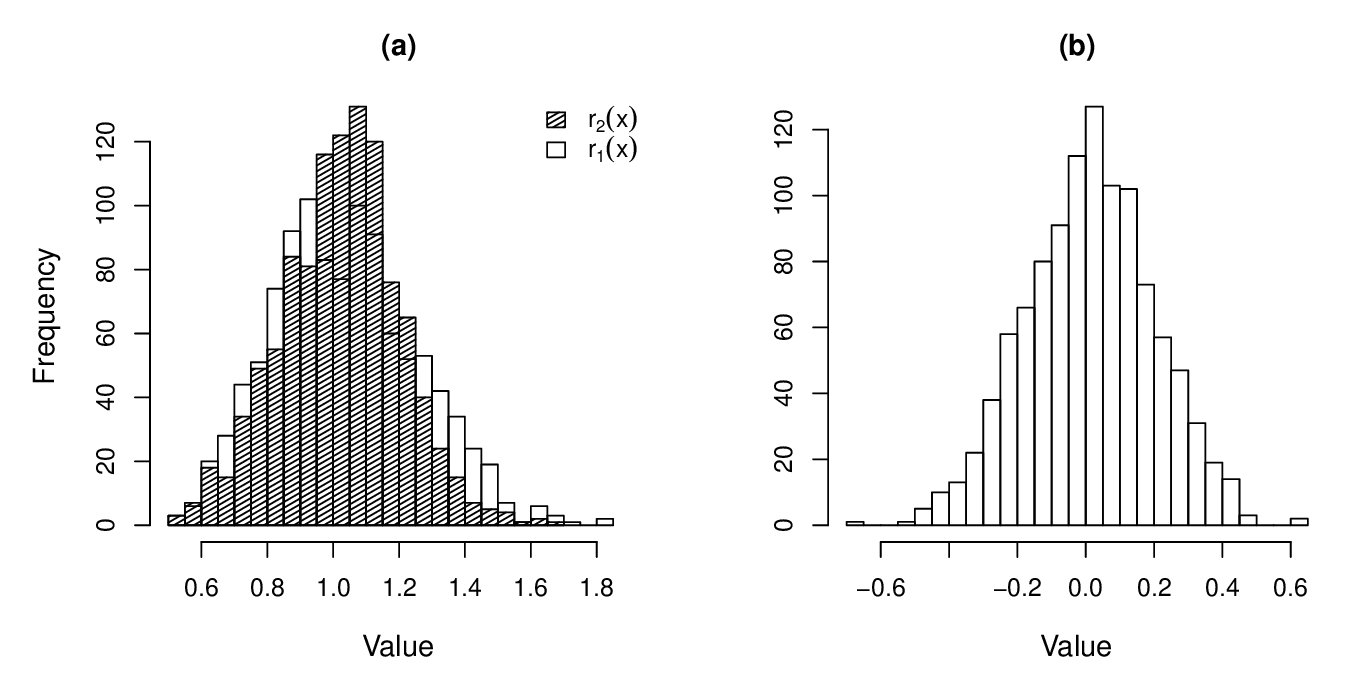}
\caption{Histograms of estimated density ratios $r_1(x)$ and $r_2(x)$ (plot (a)),
and  their difference $r_1(x)-r_2(x)$  (plot (b))
evaluated at the covariates of the test data.
}
\end{figure}

When $k=2$, the GCS  assumption  \eqref{GCS} is equivalent
to $r_1(x) - r_2(x) =  0$ for all $x$,
where $r_y(x) = Q_{X\mid Y}(x\mid y)/P_{X\mid Y}(x\mid y)$ for $y=1, 2$.
To verify the GCS assumption, we need to ensure that
$r_1(x) - r_2(x)  =  0$ for all $x$.
We randomly divide the dataset into two halves,
designating one half for training and the other for testing.
For each $y=1$ and 2, we combine the data $X_{Q}$ and $X_{P}$ with the same $Y$ value
 and estimate the density ratio
$r_y(x)$ using a logistic regression model  \citep{qin1998inferences}.
Figure \ref{fig-real-test-GCS} presents histograms depicting
the estimated density ratios $r_1(x) $ and $ r_2(x)$,
as well as their difference $r_1(x) - r_2(x)$,
evaluated at the covariates of the test data.
It is evident from Figure \ref{fig-real-test-GCS}(a)
 that   the histograms of $r_1(x) $ and $ r_2(x)$
are highly similar.
Furthermore, Figure \ref{fig-real-test-GCS}(b) reveals that
most of their  difference values are close to zero.
Both  observations lend  support to our GCS assumption.

\begin{figure}[htbp]
	\label{fig-real-result}
	\centering
	\includegraphics[width=0.8\textwidth]{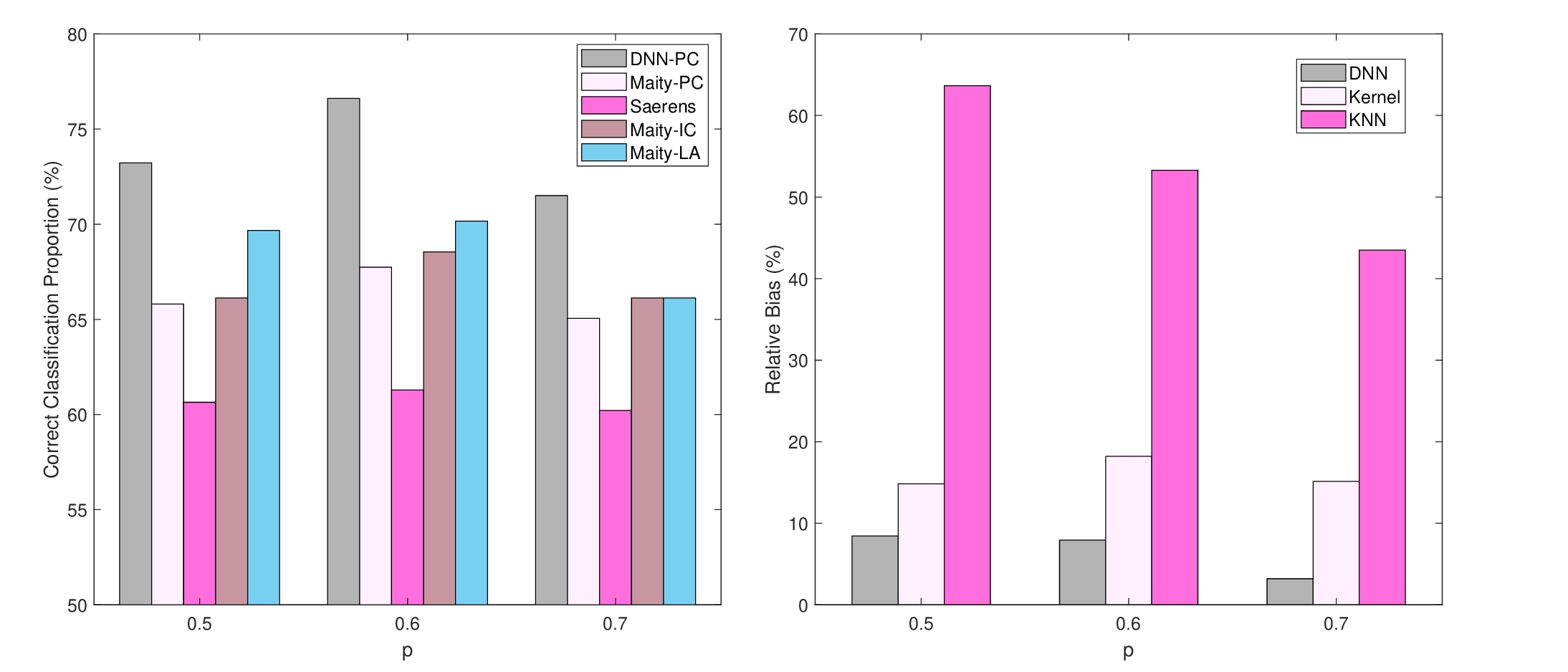}
	\caption{Plots of analysis results for the Alzheimer's Disease Dataset when $p=0.5$, $0.6$, or $0.7$.
Left plot: Correct classification proportions of the five classifiers under comparison.
Right plot: Absolute relative bias of DNN, Kernel and KNN with respect to the proportion of label $1$ in the target data. }
\end{figure}

We apply the proposed estimation and classification method and its competitors presented in Section \ref{subsection-simu-method} to predict the diagnosis results of the patients in the target domain. The samples in the target domain are randomly partitioned into training and validation sets, with the training set consisting of $100p\%$  of the total target data, where  $p=0.5$, $0.6$, or $0.7$. Although the target data here are all labeled, we consider the labels of the training set of the target data as missing when applying DNN-PC, Maity-PC, and Saerens. The correct classification proportion of each classifier is evaluated based on the validation data.
The left plot of Figure \ref{fig-real-result} shows the proportions of correct classification or correct diagnosis of the five classifiers under comparison for the three choices of $p$. It is evident that our DNN-PC consistently has the largest correct classification proportion among all five classifiers, indicating that our classifier yields the most accurate diagnosis results.
The right plot of Figure \ref{fig-real-result} presents the absolute relative biases of the three estimators (DNN, Kernel, and KNN) of $\pi_{Q,1}$ when we take the proportion of label 1 in the target data as the true value of $\pi_{Q,1}$. The plot reveals that the absolute relative bias of our DNN-based PMLE (i.e., DNN) is significantly smaller than those of the Kernel and KNN estimates. In particular, for different choices of $p$, the absolute relative biases of DNN are no greater than 10\%, whereas those of both Kernel and KNN are no less than 15\%. These results once again affirm the superiority of the proposed classifier and proportion estimation method over the existing competitors.

\section{Discussion}
\label{sec-con}

Deep learning and transfer learning have exhibited their potential across numerous domains. Under a novel GCS assumption, this paper explores the domain of deep transfer learning within the framework of classification that is based on a group of labelled source data and unlabelled target data. We meticulously establish concentration bounds for the proposed DNN estimator of the classification function and the PMLE of the class proportion of the target domain. We construct a Bayes classifier for the target domain using these two estimators and prove that it is minimax optimal in terms of excess-risk. Thanks to the remarkable approximation capability of DNN, the proposed classifier can circumvent the notorious curse-of-dimensionality issue and showcases impressive numerical performance regarding both accuracy and efficiency.

In our analysis of the Alzheimer's Disease dataset, a histogram-based visual comparison method is employed to justify the GCS assumption; however, it lacks the rigor of quantitative analysis. This thereby raises the question of how to construct a valid hypothesis test for checking the GCS assumption based on labelled source data and unlabelled target data.
Intuitively, a test statistic can be constructed by measuring the deviation from GCS. For instance, we might use the overall absolute difference between $\hat r_1(x), \ldots, \hat r_{k}(x)$, where
$ r_y(x)= Q_{X\mid Y}(x\mid y)/P_{X\mid Y}(x\mid y)$
and $\hat r_y(x)$  is its estimate. Nevertheless, it is quite challenging to derive the limiting distribution of this measure under the GCS assumption. We may leave this issue for future research.

\section*{Supplementary information}

The supplementary material  contains proofs of all the theorems and lemmas,
and  provides additional simulation results.

\section*{Acknowledgements}

Liu's research is supported by  the National Natural Science Foundation of China (12371293, 12171157,  32030063),
Fundamental Research Funds for the Central Universities and the 111 Project (B14019).

\bibliographystyle{natbib}
\bibliography{Classification-GCS-DNN}

\end{document}